\documentclass[review,hidelinks,onefignum,onetabnum]{siamart220329}

\ifpdf
\hypersetup{
  pdftitle={Constrained Local Approximate Ideal Restriction for Advection-Diffusion Problems},
  pdfauthor={Ahsan Ali, James J. Brannick, Karsten Kahl, Oliver A. Krzysik, Jacob B. Schroder, and Ben S. Southworth}
}
\fi

\ifpdf
  \DeclareGraphicsExtensions{.eps,.pdf,.png,.jpg}
\else
  \DeclareGraphicsExtensions{.eps}
\fi

\newsiamremark{remark}{Remark}
\newsiamremark{hypothesis}{Hypothesis}
\crefname{hypothesis}{Hypothesis}{Hypotheses}
\newsiamthm{claim}{Claim}

\headers{}{}

\title{Constrained Local Approximate Ideal Restriction for Advection-Diffusion Problems\thanks{Submitted to the editors June 30, 2023.
\funding{This work was funded by NSF grants DMS-2110917 and DMS-2111219.} }}

\author{Ahsan Ali\thanks{Department of Mathematics and Statistics, University of New Mexico, Albuquerque, NM, USA (\email{ahsan@unm.edu})}
\and James J. Brannick\thanks{Department of Mathematics, Penn State University, University Park, PA, USA (\email{jjb23@psu.edu})}
\and Karsten Kahl\thanks{School of Mathematics and Natural Sciences, Bergische Universit\"{a}t Wuppertal, Wuppertal, Germany (\email{kkahl@uni-wuppertal.de})}
\and Oliver A. Krzysik\thanks{Department of Applied Mathematics, University of Waterloo, Waterloo, Ontario, Canada (\email{okrzysik@uwaterloo.ca})}
\and Jacob B. Schroder\thanks{Department of Mathematics and Statistics, University of New Mexico, Albuquerque, NM, USA (\email{jbschroder@unm.edu})}
\and Ben S. Southworth\thanks{Los Alamos National Laboratory, Los Alamos, NM, USA (\email{southworth@lanl.gov}) }}

\usepackage{amsopn}

\usepackage{lipsum}
\usepackage{graphicx}
\usepackage{subcaption}
\usepackage{epstopdf}
\usepackage{algorithmic}
\usepackage[utf8]{inputenc}
\usepackage{amsmath}
\usepackage{amssymb,mathtools}
\usepackage{amsfonts}

\usepackage{cleveref}
\usepackage{booktabs}
\usepackage[total={6.25in, 9.3in}]{geometry}
\usepackage{url}

\usepackage{xcolor}

\newcommand{\aff}[1]{\ensuremath{A_{ff}^{(#1)}}} 
\newcommand{\affinv}[1]{\ensuremath{A_{ff}^{(#1),-1}}} 
\newcommand{\afc}[1]{\ensuremath{a_{fc}^{[#1]}}} 

\newcommand{\w}[1]{\ensuremath{w^{[#1]}}} 
\newcommand{\witer}[2]{\ensuremath{w^{[#1],#2}}} 
\newcommand{\witerupdate}[2]{\ensuremath{\bar{w}^{[#1]}}} 
 
\newcommand{\B}[1]{\ensuremath{B^{[#1]}}}

\newcommand{\clair}{C$\ell$AIR}
\newcommand{\RR}{\ensuremath{\mathbb{R}} } 
\newcommand{\pideal}{\ensuremath{P_{\textnormal{ideal}}}}
\newcommand{\rideal}{\ensuremath{R_{\textnormal{ideal}}}}

\graphicspath{ {./plots/} }
\DeclareMathOperator*{\argmin}{argmin}

\newcommand{\squishlist}{
   \begin{list}{--}
    { \setlength{\itemsep}{1pt}      \setlength{\parsep}{0pt}
      \setlength{\topsep}{0pt}       \setlength{\partopsep}{0pt}
      \setlength{\labelwidth}{1em}   \setlength{\leftmargin}{20pt}
      \setlength{\labelsep}{0.5em}   \setlength{\parskip}{0pt} } }

\newcommand{\squishend}{
    \end{list}  }
\nolinenumbers
\begin{document}

\maketitle 

\begin{abstract}
This paper focuses on developing a reduction-based algebraic multigrid (AMG) method that is suitable for solving general (non)symmetric linear systems and is naturally robust from pure advection to pure diffusion. Initial motivation comes from a new reduction-based algebraic multigrid approach, $\ell$AIR (local approximate ideal restriction), that was developed for solving advection-dominated problems. Though this new solver is very effective in the advection dominated regime, its performance degrades in cases where diffusion becomes dominant. This is consistent with the fact that in general, reduction-based AMG methods tend to suffer from growth in complexity and/or convergence rates as the problem size is increased, especially for diffusion dominated problems in two or three dimensions. Motivated by the success of $\ell$AIR in the advective regime, our aim in this paper is to generalize the AIR framework with the goal of improving the performance of the solver in diffusion dominated regimes. To do so, we propose a novel way to combine mode constraints as used commonly in energy minimization AMG methods with the local approximation of ideal operators used in $\ell$AIR. The resulting constrained $\ell$AIR (\clair) algorithm is able to achieve fast scalable convergence on advective and diffusive problems. In addition, it is able to achieve standard low complexity hierarchies in the diffusive regime through aggressive coarsening, something that has been previously difficult for reduction-based methods.

\end{abstract}



\section{Introduction}
We design and analyze a reduction-based algebraic multigrid (rAMG) algorithm for linear systems of algebraic equations
\begin{equation}\label{eq:mod}
A \bf{u} = \bf{f},
\end{equation} 
where $A \in \mathbb{R}^{n\times n}$
is assumed to be a sparse nonsymmetric matrix. Our focus in this paper is on solving \eqref{eq:mod} arising from discretizations of advection-diffusion-reaction partial differential equations (PDEs), which arise in various practical applications and also serve as interesting initial model problems for testing nonsymmetric AMG solvers. Ultimately, our objective is a method that is naturally robust and efficient in both the advection and diffusion limits.

Multigrid methods for solving~\eqref{eq:mod} use a relaxation process or smoother, defined in this paper by $M$, as a local solver for a sequence of coarse-level systems of equations to reduce the global error resulting from applying relaxation on the finest level.  In AMG, a recursive two-level point of view is often used, both in terms of the development of the AMG setup algorithm as well as the analysis of the solver it produces. In this two-level context, the idea is to analyze two complementary processes to efficiently solve sparse linear systems, a relaxation scheme on the fine-level, with corresponding
error propagation matrix given by $I - M^{-1}A$, and a coarse-level correction, with error propagation matrix given by $ I - P A_c^{-1} R A$. Here, $P \in \mathbb{R}^{n \times n_c}$, $n_c < n$, denotes the interpolation matrix that maps corrections from the coarse-level, $R \in \mathbb{R}^{n_c \times n}$ is the restriction matrix that maps residuals to the coarse-level, and $A_c = RAP$ is the coarse-level operator. The error propagation matrix of the resulting two-level method, from which a multilevel method is defined recursively, reads
\begin{equation}\label{eq:tl}
E_{TG} = (I-M^{-1}A)(I-P A_c^{-1} R A).
\end{equation} 
The AMG solver is then defined from this two-level scheme by recursively applying it on the coarse-level to approximate $A_c^{-1}$, and reduce the errors remaining after applying relaxation on the finer levels.

In AMG, the smoother $M$ is typically fixed to be a simple point-wise method and then $R$ and $P$ are constructed in
an automated setup algorithm that takes as input the system matrix $A$.  In this paper, since we consider a reduction-based AMG framework, we assume that the smoother is a point-wise $F$-relaxation scheme, where the set $F$ denotes the set of fine variables in the coarse-fine ($C/F$) splitting of the degrees of freedom $\Omega$ such that $C\cup F = \Omega$ and $C\cap F = \emptyset$.  
In this setting, the main tasks in the AMG setup algorithm are to construct the restriction and interpolation matrices $R$ and $P$ such 
that certain {\em approximation properties} hold and $R$ and $P$ (and thus $A_c=RAP$) are {\em sparse} matrices.  The latter sparsity requirement implies 
that the setup procedure can be efficiently applied recursively to $A_c = RAP$ in order to construct an optimal multilevel AMG
solver.

In the SPD case, the variational choice $R=P^T$ is used and the weak approximation property for $P$ bounds the convergence rate of the two-level method in the $A$-norm.  Notably, the weak approximation property has as its minimizer the so-called ideal interpolation
matrix~\cite{Falgout}.  This ideal form of interpolation gives rise to the Schur-complement of $A$ as
the coarse-level system matrix and, thus, coincides with block Gaussian elimination based on the given coarse-fine splitting.
In contrast to classical AMG methods that use a global smoother in the AMG solver (e.g., weighted-Jacobi or lexicographic Gauss-Seidel), reduction based methods that are motivated by this block factorization interpretation choose smoothers that focus on the subspace defined by the fine variables $F$. 

Traditionally, classical AMG methods are very effective for sparse SPD problems, e.g., discretizations of various diffusion and heat conduction problems, and not as effective for nonsymmetric problems, whereas, reduction-based AMG methods have been developed that work well for nonsymmetric problems with a near-triangular (upwinded) structure, e.g., space-time discretizations \cite{st-air} and advection-dominated PDEs \cite{Ben,air1}. The parallel-in-time method called multigrid reduction in time (MGRIT) \cite{mgrit} is another recent application of reduction-based multigrid methods, where time integration is recognized as suitable for reduction-based approaches because time is one-dimensional and optimal coarsening is practical; however, MGRIT differs from rAMG methods in using a non-Galerkin coarse grid typically based on rediscretization in time. Interestingly, theory \cite{mgrit-theory-ben,mgrit-theory-jacob,southworth2021tight} (and practice) has indicated that without special care (e.g., \cite{mgrit-adv}) MGRIT is effective on space-time PDEs with parabolic problems, particularly SPD spatial discretizations, but struggles with hyperbolic PDEs and highly nonsymmetric spatial discretizations. 

Indeed, solving more general non-symmetric systems arising from higher dimensional spatial PDEs or space-time PDEs is more complicated in terms of balancing the convergence and complexity of a multilevel solver.  A viable approach for solving the latter higher dimensional problems is given by the $\ell$AIR solver~\cite{Ben}. The approach has been extensively studied and tested for analogous advection-diffusion-reaction model problems we consider here~\cite{BenTheory,air1,Ben}, and also demonstrated on more complex advection-dominated physics, e.g., \cite{southworth2022fast,dargaville2023air}. In general, $\ell$AIR performs very well for advection-dominated problems, with performance more mixed for diffusive problems, i.e., for problems with strong diffusion $\ell$AIR suffers from degraded convergence and similar complexity issues as other reduction-based AMG methods for diffusion type problems in higher dimensions. We mention the recent paper~\cite{Scott} which represents the state-of-the-art in research on techniques for improving the performance of reduction-based AMG solvers for diffusion problems. There, the traditional approach of constructing a row-wise approximation of ideal interpolation is considered and the focus is on improving approximations to $A_{ff}^{-1}$ using sparse approximate inverse techniques, similar to the sparse Krylov approximations that have been studied, and very recently used for a highly efficient 
parallel variation of $\ell$AIR \cite{dargaville2023air}.

In this paper, we combine the reduction-based principles of $\ell$AIR with the mode constraints of energy-minimization AMG \cite{Mandel98energyoptimization,Wan,Brannick_Trace_06,OlScTu2011,Jacob,root-node}. The $\ell$AIR algorithm offers sparse, accurate approximations of ideal transfer operators in the strongly advective regime, which also yields excellent AMG approximation properties \cite{BenTheory}. In contrast, in diffusive regimes, accurate and sparse approximations of ideal transfer operators are generally not viable due to the density of $A_{ff}^{-1}$; because no other information is taken into account, the AMG approximation properties of $\ell$AIR also suffer or the complexity must dramatically increase. Classical and energy-minimization AMG methods suffer from the opposite problem -- they tend to offer excellent AMG approximation properties for diffusive problems, but very poor approximation properties in the advective regime \cite{BenTheory}. On a high level, this is because the near nullspace of advective operators cannot be represented simply by smoothed constant vectors, the basis for almost all classical AMG and energy minimization methods. Energy minimization methods inadvertently further block their potential by using some form of normal equations to perform energy minimization in the nonsymmetric setting (wherein the nonsymmetric matrix $A$ does not define a natural energy or minimization). Such minimization converges to the ideal transfer operators \emph{of the normal equations} \cite{root-node}, even though a sparse accurate approximation of ideal transfer operators is often viable for matrix $A$, which $\ell$AIR successfully targets.

Here, we recognize this subtlety in approximation properties that guides the success of different AMG methods in different regimes, and define a new constrained variation of the $\ell$AIR algorithm. Constrained $\ell$AIR (\clair) directly approximates the ideal transfer operators using a similar objective as in $\ell$AIR, while constraining the range of $P$ or $R^T$ to include known or expected near-null space mode(s), thereby harnessing the power of reduction-based and energy-minimization methods in their respective regimes.  In relaxation, we restrict ourselves to simple $F$- and $C$-point Jacobi relaxation, demonstrating that with careful construction of transfer operators, we are able to apply an efficient reduction-based method for diffusive-dominated problems. 
As it turns out, the construction of sparsity patterns for transfer operators in $\ell$AIR and root-node based energy minimization \cite{root-node} are very similar in principle, namely that they are defined column-wise for $P$ and row-wise for $R$. As a result, we are also able to naturally incorporate an aggressive root-node approach to choosing coarse variables~\cite{Jacob,root-node} for diffusion dominated problems, ameliorating the (necessarily) high complexity that tends to arise in reduction-based and advective solvers \cite{Ben}.

In the case of anisotropic diffusion, we observe that the energy minimization approach depends crucially on the aggressive root-node coarsening technique \cite{Jacob}. In this case, the aggressive coarsening selects a relatively small number of coarse variables to compensate for the complexity of the long-range interpolation required for the non-grid aligned anisotropy 
 \cite{Jacob}. The outcome is a low-complexity, effective AMG solver. By combining neighborhoods of fine degrees of freedom (DOFs) into a single coarse variable, root-node accomplishes this by greedy aggregation \cite{VaMaBr1996}, in which the seed point of each aggregate becomes the $C$-point (root-node) and the remaining degrees of freedom become $F$-points. The member points ($F$ and $C$) of each aggregate define the initial nonzero pattern of the corresponding column of $P$. Generally speaking, this greedy aggregation procedure places $C$-points farther apart than the traditional Ruge-St\"{u}ben (RS) $C/F$-splitting of the fine degrees of freedom. This results in a more aggressive coarsening and fewer coarse variables (i.e., lower complexity) for root-node \cite{Jacob, root-node}. The column-wise viewpoint that we employ in the definition of \clair\, interpolation allows us to similarly use aggressive root-node coarsening to control complexity, particularly for diffusive problems, and to also incorporate constraint vectors for improved AMG convergence that are fit into $\mbox{span}(P)$ (similar to  smoothed aggregation methods \cite{VaMaBr1996, root-node}).

This paper is organized as follows. 
The next section reviews reduction-based AMG and the $\ell$AIR framework that motivates our new method, and provides practical connections to energy minimization AMG methods that facilitate our method design. Our new method \clair\  is presented in \Cref{sec:iter-air}, and \Cref{sec:num} contains numerical results that illustrate the performance of our algorithms applied to various discretizations of our model advection-diffusion problem. In particular, \clair\ is able to maintain fast robust convergence on advection-dominated problems, while also yielding low-complexity scalable solutions for diffusion dominated problems. To complement the numerical results for the proposed method \clair, we present a study of classical AMG weak and strong approximation properties in Appendix A.

\section[Reduction based AMG and lAIR interpolation]{Reduction based AMG and $\ell$AIR interpolation}
Let $A\in\mathbb{R}^{n\times n}$ and assume that the degrees of freedom $\Omega = \{1,...,n\}$ are partitioned in the classical sense such that we have $n_c$ $C$-points and $n_f$ $F$-points.  Then, $A$ can be represented in the following block form:
\begin{align}
\label{eqn:Asplit}
A & = \begin{pmatrix} A_{ff} & A_{fc} \\ A_{cf} & A_{cc}\end{pmatrix}.
\end{align}
As before, define $P:\mathbb{R}^{n_c}\mapsto\mathbb{R}^n$ and $R:\mathbb{R}^{n}\mapsto\mathbb{R}^{n_c}$ as interpolation and restriction respectively. Further, assume that $C$-points are interpolated and restricted by injection in the classical AMG sense; then, the transfer operators $P$ and $R$ in  reduction based AMG can be written in the following block form:
\begin{align}
\label{eqn:basicsplit}
P = \begin{pmatrix} W \\ I\end{pmatrix},
\hspace{3ex} R = \begin{pmatrix} Z & I \end{pmatrix},
\end{align}
where the ordering here is useful in formalizing a reduction-based AMG method.
Design of classical reduction based AMG methods is motivated by the observation that the ideal interpolation operator is the unique operator 
\begin{align}
\pideal & = \begin{pmatrix} -A_{ff}^{-1}A_{fc} \\ I \end{pmatrix}\label{eq:p_ideal}
\end{align}
that eliminates the contribution of the coarse-grid correction $\mathbf{e}_c$ to the $F$-point residual:
\begin{align}
AP\mathbf{e}_c = \begin{pmatrix}\mathbf{0}\\ S\mathbf{e}_c \end{pmatrix} \label{eq:def1}. 
\end{align}
Assuming $R$ and $P$ take the form of \eqref{eqn:basicsplit}, we have independent of $A$ that the Petrov-Galerkin coarse grid satisfies $RAP = S\coloneqq A_{cc}-A_{cf}A_{ff}^{-1}A_{fc} $, where $S$ is the Schur complement \cite{Falgout,Ben}.  

Reduction based methods for diffusion type problems have traditionally assumed the classical AMG form of interpolation in~\eqref{eq:p_ideal} and approximate $-A_{ff}^{-1}A_{fc}$ by solving for each $i\in C$
\begin{align}
A_{ff} W^{[i]} = -A_{fc}^{[i]}, 
\end{align}
for the $n_c$ columns of the interpolation weight matrix 
$W$ defined as in~\eqref{eqn:basicsplit}.
Here, we use the notation $X^{\{i\}}$ to refer to the $i$th row of a matrix $X$, whereas, $X^{[i]}$ will denote a column of the matrix.
In this classical reduction-based setting, choosing the coarse grid degrees of freedom $C$ 
as well as choosing the sparsity structure of the rows of $W$ are done using classical AMG coarsening and strength of connection heuristics and the resulting algorithms typically lead to high grid and operator complexities, even when very simple approximations of $A_{ff}^{-1}$ are used.  
We mention that various approximations to $A_{ff}^{-1}$ in the computation of $P$ and in relaxation are possible and have been considered in the literature, and we reiterate that in our development we focus on designing an approach that achieves both low grid and operator complexities while at the same time only requiring the simplest diagonal (Jacobi) F-relaxation scheme for fast convergence.

\subsection[Review of lAIR]{Review of $\ell$AIR}
The $\ell$AIR approach that we build our new method around is based on the similar observation that the ideal restriction operator is the unique operator \cite{Ben}
\begin{align}
\rideal & = \begin{pmatrix} -A_{cf}A_{ff}^{-1} & I \end{pmatrix}\label{eq:r_ideal}
\end{align}
that eliminates all error at $F$-points:
\begin{align}
RA\begin{pmatrix}\delta\mathbf{e}_f\\\mathbf{0}\end{pmatrix} = \mathbf{0}\hspace{3ex}\forall \delta\mathbf{e}_f. \label{eq:def0}
\end{align}
Here, the ordering of the equations is again based on a splitting of $\Omega = C \cup F$. 
The $\ell$AIR  approach is then based on setting $RA$ equal to zero in~\eqref{eq:def0} within a pre-determined $F$-point sparsity pattern for $Z$. A similar method to approximate ideal interpolation can be expressed as satisfying $AP=0$ exactly a specified sparsity pattern for $W$.  Note that the AIR approach can also be seen as directly approximating the action of $\rideal$ on F-points, where $\rideal A = (\mathbf{0} , S)$, for the Schur complement $S$. Expressing this result in terms of some matrix $Z$, the approach is equivalent to satisfying
\begin{equation}\label{eq:z-formula}
    ZA_{ff} = -A_{cf}
\end{equation}
within a predetermined sparsity pattern for $Z$. Here, the AIR approach is clearly different from the classical (reduction) AMG \cite{1985BrandtA_McCormickS_RugeJ-aa,JWRuge_KStuben_1987a,maclachlan2006adaptive,FromKahlMaclZikaBran2010} approach in that~\eqref{eq:z-formula} involves solving for the $n_c$ rows of $Z$ (or $R$), which now gives a column-wise view of computing $R^T$ and thereby an $\ell$AIR style form of $P$.~\footnote{We note that the reduction-based AMG approaches \cite{maclachlan2006adaptive,FromKahlMaclZikaBran2010} have also previously explored the use of constraint vector(s), however in addition to our different column-wise view of computing $R^T$, these approaches use cheap approximations to $A_{ff}^{-1}$ (e.g., diagonal) with a more expensive multilevel adaptive approach for generating the constraint vector.  In this work, we take the more expensive $\ell$AIR approach to approximating $A_{ff}^{-1}$ based on small block inverses, but combine that with a cheap constraint vector similar to energy-minimization methods, which require only a few relaxation sweeps on each level for improvement.}

Denoting indices of the sparsity pattern for the $i$th row of $Z$ as $\mathcal{Z}_i = \{\ell_1,...,\ell_{S_i}\}$, where $S_i = |\mathcal{Z}_i|$ is the size of the sparsity pattern, the resulting (transposed) linear system for $Z$ takes the form
\begin{align}
\begin{pmatrix} a_{\ell_1 \ell_1} & a_{\ell_2\ell_1} & ... & a_{\ell_{S_i}\ell_1} \\
a_{\ell_1\ell_2} & a_{\ell_2\ell_2} & ... & a_{\ell_{S_i}\ell_2} \\
\vdots & & \ddots & \vdots \\
a_{\ell_1\ell_{S_i}} & a_{\ell_2\ell_{S_i}} & ... & a_{\ell_{S_i}\ell_{S_i}}\end{pmatrix}
\begin{pmatrix} z_{i\ell_1} \\ z_{i \ell_2} \\ \vdots \\ z_{i\ell_{S_i}}\end{pmatrix}
& = - \begin{pmatrix} a_{i\ell_1} \\ a_{i\ell_2} \\ \vdots \\ a_{i\ell_{S_i}}\end{pmatrix} \label{eq:system}.
\end{align}
As demonstrated in \cite{Ben,air1}, for upwinded advection-dominated problems, very good sparse approximations to $A_{ff}^{-1}$ can be made, which means that satisfying \eqref{eq:z-formula} within a sparsity pattern provides an accurate approximation to $\rideal$. For advection-dominated problems, this approach also provides very good approximation properties of the resulting restriction operator, as demonstrated in ~\cite{BenTheory} and the numerical results. In contrast, relying only on solving \eqref{eq:z-formula} for a diffusion dominated problem is unlikely to be effective, because $A_{ff}^{-1}$ is generally more dense in this setting and not well approximated by a sparse matrix. Hand-in-hand with this, the resulting approximation properties are also poor. This has been mitigated reasonably well by combining $\ell$AIR restriction with classical AMG interpolation for problems with strong diffusion, but the complexity remains high and convergence sub-par compared with classical AMG or energy minimization methods. 

\subsection{Review of Root-node AMG}
\label{sec:rn_amg}
One way to conceptualize root-node AMG \cite{root-node} is as a combination of classical and aggregation-based multigrid methods. Root-node AMG employs a hybrid strategy in which smoothed aggregation type strength-of-connection is used and aggregates are created using standard aggregation routines. One node is selected as the `root-node' in each aggregate which corresponds to a $C$-point, while all other nodes in the aggregate are identified as $F$-points. After that, transfer operators are created in the following manner. 

Root-node utilizes algebraically smoothed candidates $B$, which are fit into the span of interpolation, and in order to recurse, coarse versions, $B_c$, of the candidates are obtained by injecting $B$ to the $C$-points. Next, the initial tentative interpolation $T$ is formed by injecting only the first $q$ candidates over each aggregate, where $q$ is the block size of the original matrix ($q=1$ for a scalar problem). On the coarse grid, each root-node thus represents $q$ DOFs.
On $T$, a further step is taken to normalize every column. This procedure produces the following form 
\begin{equation}\label{eq:root-nodeform}
  T =
  \begin{array}{c@{}c}
    
      \begin{pmatrix}
        W \\
        I\\
      \end{pmatrix}
   
  &
  \begin{array}{l}
    \} \,\,\textnormal{Non Root-nodes}\\
    \} \,\,\textnormal{Root-nodes}\\
  \end{array}
  \end{array}
\end{equation}
For $q=1$, $T$ has non-overlapping columns; for $q>1$, $W$ is block diagonal. 

The remaining candidates are projected into $\mbox{range}(T)$ in the Euclidean inner-product if there are more than $q$ candidates. It is expected in root-node AMG that the allowed sparsity pattern of $T$ has enough DOFs to make this an underdetermined problem.  The sparsity pattern of $T$ is typically grown with strength of connection information, as we later do for \clair.
Consequently, a minimal norm update is applied to each row of $T$, guaranteeing that $T B_{c} = B$ and $T$ adheres to the sparsity pattern.  The interpolation $P$ is then generated using subsequent energy-minimization updates to $T$ (briefly covered in the following subsection). Root-node minimizes energy by solving $A\begin{pmatrix}
    W \\
    I
\end{pmatrix}=0$ subject to the mode interpolation constraints so that the solution is non zero.

\subsection{Energy minimization and mode constraints}
\label{sec:enmin_modeconst}

Another well-known class of AMG methods is that of energy-minimization. We will particularly focus on root-node energy-minimization \cite{root-node}, which shares the $C/F$-splitting design of reduction-based methods. One key part of energy-minimization is the use of constraints during the minimization process.  For efficiency reasons, the sparsity pattern of $P$ is constrained.  Let $\mathcal{W}$ be the sparsity pattern for the F-rows in $P$.  We denote that $P$ obeys the sparsity pattern constraint with 
\begin{equation}
\label{eqn:sparsityconst}
P \in \begin{bmatrix} \mathcal{W}\; I \end{bmatrix} \quad \mbox{or} \quad W \in \mathcal{W}.
\end{equation}
For approximation property reasons, a near nullspace mode constraint is also typically enforced where
\begin{equation}\label{eq:Pcontraint}
    B \in \mbox{span}(P),
\end{equation}
and $B \in \RR^{n,k}$ is a set of $k$ global near nullspace modes.  With $P$ of the form~\eqref{eqn:basicsplit} and equation~\eqref{eq:Pcontraint}, this then implies that $B_c = [0\; I]\, B$, i.e., the fine-grid $B$ is injected to the coarse-grid, along with the constraint,  
\begin{equation}
   \label{eqn:global_const1}
    P B_c = B.
\end{equation}

Similar to $\ell$AIR, energy-minimization AMG takes a column-oriented view. Letting $P^{[j]}$ denote column $j$, construction of transfer operators is based around a minimization along the lines of
\begin{equation}
\label{eqn:rn_min}
P = \argmin_P \sum_j \| P^{[j]} \|_{\mathcal{X}}^2 \mbox{, such that constraints (\ref{eqn:global_const1}) and (\ref{eqn:sparsityconst}) are satisfied,}
\end{equation}
where $\mathcal{X}$ denotes some norm, usually $A$ for SPD operators or $A^*A$ for nonsymmetric matrices. This column-oriented view is analogous to $\ell$AIR in \eqref{eq:system}, but here we are minimizing in some energy-induced inner product, rather than solving each block equation exactly as in $\ell$AIR. 

Interestingly, this distinction has more profound consequences for the efficacy of the methods. For SPD operators, root-node uses projected conjugate gradients (CG) for equation (\ref{eqn:rn_min}); without constraints (\ref{eqn:global_const1}) and (\ref{eqn:sparsityconst}), such a minimization procedure converges to \pideal\ for $A$ \cite[Lemma 4.2]{root-node}. However in the nonsymmetric case, root-node uses a projected generalized minimal residual (GMRES) for equation \eqref{eqn:rn_min}, which in turn (without constraints (\ref{eqn:global_const1}) and (\ref{eqn:sparsityconst})), converges to $\pideal$ for $A^* A$ \cite[Lemma 4.6]{root-node}. Indeed, by posing the unconstrained problem as an (overdetermined) minimization, root-node is required to define an energy-induced inner product through the normal equations, which in turn leads to the approximation of $\pideal$ for the normal equations rather than directly for the operator of interest. In contrast, the base algorithm of $\ell$AIR directly approximates (and converges to) $\pideal$ for the original operator $A$. This is possible because the algorithm is built around local matrix approximation \eqref{eq:z-formula} rather than a matrix-induced norm, which does not naturally exist for non-SPD matrices. 

As it turns out, which ideal operators we are trying to approximate (without constraints) is an important distinction for highly advective problems, which typically generate discretization matrices that are close to block lower-triangular in some ordering \cite{air1}. For such cases \cite{air1}, $\ell$AIR achieves good sparse approximations to $A_{ff}^{-1}$ for computing $\rideal$ and overall excellent AMG convergence. However if one were to approximate ideal transfer operators based on $A^*A$ instead of $A$, the block lower-triangular structure that is key to achieving sparse and accurate approximations to $A_{ff}^{-1}$ and $\rideal$ is completely lost. The natural result is that a good sparse approximation to $A_{ff}^{-1}$ becomes more difficult to compute
(see \cite{air1}, Section 4 for more discussion), and the resulting AMG method is significantly less effective. Thus although the underlying problem that energy-minimization is based around, namely approximating $AP = \mathbf{0}$, is almost equivalent to $\ell$AIR, by formulating via energy minimization the resulting class of methods yield lackluster performance on highly advective problems. 

\subsection{The best of both worlds}

Looking carefully at the $\ell$AIR and energy-minimization approaches leads us to consider a new interpretation combining the best of both worlds. By directly approximating $\rideal$ of the original operator $A$, $\ell$AIR is able to construct highly effective transfer operators for advection-dominated problems; in contrast, mode constraints are fundamental to the efficacy of energy-minimization methods for diffusion dominated problems (indeed, without constraints energy minimization alone is generally not effective). Although $\ell$AIR is based around approximation of ideal transfer operators, we can also think in terms of approximation properties -- consider each row of $R$ as the local fine-grid mode being restricted to a given C-point, where these modes should be local representations of the smooth error. This is exactly what happens in classical smoothed aggregation \cite{VaMaBr1996}, as well as when bilinear interpolation is used in geometric MG for diffusion. Thus we propose a new constrained $\ell$AIR (\clair) method that is built around directly approximating the ideal transfer operators of $A$ in an $\ell$AIR framework, regardless of whether $A$ is SPD or nonsymmetric, while also incorporating mode constraints as in energy minimization to improve robustness in diffusion dominated problems. 

In summary, the proposed \clair\ approach combines strengths of root-node and $\ell$AIR, with the goal of a robust solver in both the advective and diffusive regimes.  The \clair\ approach for constructing $\ell$AIR interpolation with constraints is directly related to solving \eqref{eq:z-formula} and thus is a constrained approximation of \pideal\ for the original operator, in contrast to root-node, which targets the normal equations in the nonsymmetric setting. Also, all of our proposed generalizations can be used to build $R$ and/or $P$.

\section[]{Constrained $\ell$AIR transfer operators}\label{sec:iter-air}

The main goals of our new reduction-based AMG method built around $\ell$AIR-style interpolation, is to have a solver that (i) works well for both advection and diffusion problems, (ii) allows for the incorporation of mode interpolation constraints (local or global), and (iii) controls complexity in a reduction setting for diffusive problems through an aggressive root-node coarsening. 
To achieve these goals, we consider mixing ideas from $\ell$AIR, which works well for advection, with energy-minimization and smoothed aggregation (SA) \cite{VaMaBr1996}, which work well for anisotropic diffusion problems and allow for mode constraints.  Another key component of the algorithm for diffusion dominated problems, described in detail in a latter section, is our use of a root-node aggregation-based coarsening algorithm.

The overall $\ell$AIR interpolation scheme we consider is described as follows.  Similar to the energy-minimization discussion above, we will 
enforce that $P$ obeys the sparsity pattern constraint (\ref{eqn:sparsityconst}).
Regarding the sparsity pattern $\mathcal{W}$ for the $F$-rows, let $\mathcal{W}_i$ denote the sparsity pattern for the $i$th column of $W$, that is $\mathcal{W}_i = \{m_1, \dots, m_{T_i} \}$ and the number of nonzeros in column $i$ equals $T_i = |\mathcal{W}_i|$.
Define \aff{i} as $A_{ff}$ restricted (in rows and columns) to the sparsity pattern $\mathcal{W}_i$, and 
 \afc{i} and \w{i} as $A_{fc}$ and $W$ restricted to column $i$, respectively, with rows restricted to the sparsity pattern of $\mathcal{W}_i$.
The standard $\ell$AIR approach for finding each \w{i} is then equivalent to solving: for each $i\in n_c$ solve
\begin{equation}
    \label{eqn:airwi}
    \aff{i} \w{i} = - \afc{i}.
\end{equation}
If we consider solving \eqref{eqn:airwi} at each $i$, we can rewrite the procedure as the following global block diagonal system, where we assume $n_c$ points on the coarse grid:
\begin{equation}
    \label{eqn:Affdefn}
    \aff{*} \w{*} := 
    \begin{bmatrix}
    \aff{0} &           &           & \\
            & \aff{1}   &           & \\
            &           & \ddots    & \\
            &           &           & \aff{n_c}
    \end{bmatrix}
    \begin{bmatrix}
    \w{0} \\ \w{1} \\ \vdots \\ \w{n_c} 
    \end{bmatrix}
    =
    -
    \begin{bmatrix}
    \afc{0} \\ \afc{1} \\ \vdots \\ \afc{n_c}
    \end{bmatrix}
    = - \afc{*}.
\end{equation}
The solution of system \eqref{eqn:Affdefn} gives us classic AIR interpolation. As a side note, let $W$ be stored in a sparse matrix format, then the global vector of all the nonzeros of $W$, denoted \w{*}, corresponds to the data array for the sparse representation of $W$ when stored in compressed column format.

\subsection{Proposed Method with Global Constraints}
Given the above formulation of $\ell$AIR in \eqref{eqn:Affdefn}, we define the procedure for incorporating a global mode interpolation constraint into the approach. To enforce the mode interpolation constraint (\ref{eqn:global_const1}), define the matrix $Q$, such that $Q \w{*}$ is equivalent to $P B_c$, i.e., the entries of $B_c$ populate $Q$ such that the constraint equation \eqref{eqn:global_const1} is equivalent to saying
\begin{equation}
    \label{eqn:Qdef}
    Q \w{*} = \begin{bmatrix} \B{0}|_F \\ \B{1}|_F \\ \vdots\\ \B{k}|_F \end{bmatrix} = \B{*}|_F,
\end{equation}
where $\B{i}$ is the $i$th column of $B$, $\B{*}$ represents the columns of $B$ stacked vertically, and  $\B{*}|_F$ represents columns of $B$ stacked vertically, but restricted to the fine-grid points $F$.  Equation \eqref{eqn:Qdef} is also equivalent to the constraint that $W B_c = B|_F$.

Thus our constrained minimization problem is
\begin{subequations}
\begin{align}
\min_{\w{*}} \left\| \afc{*}  + \aff{*} \w{*} \right\|_2  \label{eqn:minproba} \\
\mbox{subject to } Q \w{*} = \B{*}|_F.  \label{eqn:minprobb}
\end{align}
\end{subequations}
This is an equality constrained minimization problem, with various 
solution approaches \cite{GoHrNo2001}.

\subsubsection{Direct Solution to Minimization Problem}
The minimization problem (\ref{eqn:minproba})--(\ref{eqn:minprobb}) can be solved directly via
the following Karush–Kuhn–Tucker (KKT) system
\begin{equation}
    \label{eqn:KKT}
    \begin{bmatrix}
    (\aff{*})^T \aff{*} & Q^T \\
    Q       &  \mathbf{0}   \\
    \end{bmatrix}
    \begin{bmatrix}
    \w{*} \\ \lambda
    \end{bmatrix}
    =
    \begin{bmatrix}
    - (\aff{*})^T \afc{*} \\ \B{*}|_F
    \end{bmatrix}.
\end{equation}
The (1,1) block of equation (\ref{eqn:KKT}) uses the normal equations, as the minimization principle requires
an SPD matrix and we do not assume that $\aff{*}$ is SPD.

System \eqref{eqn:KKT} could be solved exactly via a Schur complement approach.  If this is done, letting $\bar{A}^{-1} = \left( (\aff{*})^T \aff{*} \right)^{-1}$, the solution is
\begin{align}
    \label{eqn:KKTsoln}
    \w{*} = & \left(I - \bar{A}^{-1} Q^T( Q \bar{A}^{-1} Q^T)^{-1} Q\right) \bar{A}^{-1} (\aff{*})^T \afc{*}\; + \; \left( \bar{A}^{-1} Q^T (Q \bar{A}^{-1} Q^T)^{-1} Q^T \right)^{-1} \B{*}|_F,
\end{align}
where $Q$ is the rectangular constraint matrix from above. Upon inspection, computing (\ref{eqn:KKTsoln}) is an expensive endeavour, especially $( Q \bar{A}^{-1} Q^T)^{-1}$ and potentially, the transpose of $\aff{*}$. 
To avoid these costs, we consider an iterative approach related to energy-minimization AMG \cite{OlScTu2011}. 

\subsubsection{Iterative Solution to Minimization problem}

A review of approaches for obtaining an inexpensive iterative solution to~\eqref{eqn:KKT} is given in~\cite{GoHrNo2001}.  One option previously used for AMG (e.g., for root-node) is projected Krylov methods. Here, an initial guess (tentative interpolation) $\witer{*}{t}$ that satisfies the constraints is constructed,
so that $Q \witer{*}{t} = \B{*}|_F$. Then, a projected Krylov method using $Q$ is applied to solve the interpolation equation (\ref{eqn:Affdefn}). The inverse $( Q \bar{A}^{-1} Q^T)^{-1}$ is not required and the transpose is not needed, because we only compute the residual for the interpolation equation (\ref{eqn:Affdefn}) when computing a descent direction for equation (\ref{eqn:minproba}).  The previous works \cite{OlScTu2011,root-node} use such a projected CG and GMRES approach for the symmetric and nonsymmetric cases, respectively.  However as noted in Section \ref{sec:enmin_modeconst}, the nonsymmetric GMRES approach will approximate \pideal\ in the constraint space for the normal equations, which is not desirable.

Thus here, we consider a simpler and cheaper linear
iteration for minimizing (\ref{eqn:minproba}) that approximates \pideal\ for the original operator in the constraint space.  We find that this approach yields effective restriction and interpolation operators.  

An additional cost consideration is whether or not to precondition such an iterative solve.  Since the 
matrix $\aff{*}$ is block diagonal, each block could be inverted (or approximately inverted). Thus,
we consider the use of approximate inverse preconditioners of the following form:
\begin{equation}
    \affinv{*} \approx \widehat{\affinv{*}} =
    \begin{bmatrix}
    \widehat{\affinv{0}} &                        &           & \\
                         & \widehat{\affinv{1}}   &           & \\
                         &                        & \ddots    & \\
                         &                        &           & \widehat{\affinv{n_c}}
    \end{bmatrix},
\end{equation}
where $\widehat{\affinv{i}}$ represents an approximate inverse to that block.  Importantly, this inverse is \textit{local} and can be generated in a variety of ways, e.g., diagonal, GMRES, or ILU approximations to each individual block inverse.  This is \textit{in contrast} to the ``classic" energy-minimization which uses a single global Krylov polynomial to simultaneously approximate all block inverses. It is our belief that the block inverse approach is more effective. In particular, a global Krylov polynomial effectively assumes that each block has the same minimizing polynomial, whereas in reality each block will likely have its own unique minimizing polynomial (distinct from other blocks). The local approach allows us to calculate more accurate local inverses faster through locally accurate approximations and polynomials, or to solve each local equation directly.

\subsubsection[]{Proposed Algorithm for Computing $R$ and $P$}
\label{sec:prop_algorithm}

We now present our simple iterative scheme for minimizing
(\ref{eqn:minproba})--(\ref{eqn:minprobb}) in Algorithm \ref{alg:iterativeAIR}.
The approach is a projected one-step iteration, which iteratively finds
AIR-like interpolation operators with constraints. For restriction, the simplest approach on paper applies Algorithm \ref{alg:iterativeAIR} to $A^T$. However if forming a transpose is computationally expensive, one can also reformulate Algorithm \ref{alg:iterativeAIR} relative to $R A = 0$, as in the original $\ell$AIR method \cite{Ben}. Then, the algorithm will still extract small submatrices $\aff{i}$, but after the extraction these submatrices will be transposed. 

As input, the algorithm takes the operator $A$ and corresponding strength of
connection matrix $S$.  Unless noted otherwise, we use the classical strength
measure \cite{JWRuge_KStuben_1987a}.  The input sparsity degree pattern $m$
determines how wide the sparsity pattern in $P$ will be, with $m=1$
corresponding to distance-one interpolation based on the sparsity pattern of
$S_{fc}$.  Most commonly, we will use $m=2$, which expands the sparsity to
consider degree-two connections, similarly to $\ell$AIR and root-node.  Next,
the input ``Coarsen type" considers whether a classical ``FC" coarsening, e.g.,
Ruge-St\"{u}ben coarsening \cite{JWRuge_KStuben_1987a}, or an aggregation-based
coarsening is used.

The coarsen type controls the base sparsity pattern for $P$.  If classical FC
coarsening is used, then the base pattern $T$ comes from the FC rows of the
strength matrix, $S_{fc}$. If an aggregation-based coarsening is used, which is
significantly more aggressive, then the ``Aggregation Operator"\footnote{\label{foot:agg} The
``Aggregation Operator" is generated by the algorithm from \cite{VaMaBr1996}.
First, an aggregation (disjoint splitting) is computed with a greedy graph
algorithm that finds the next degree-of-freedom (root node) with all unmarked
neighbors and places those degrees-of-freedom in an aggregate and then marks
them. A clean-up phase takes all unmarked degrees-of-freedom and places them in
an adjacent aggregate.  The root node of each aggregate is treated as a
C-point.  This procedure produces a $C/F$-splitting that is significantly more
aggressive (fewer C-points) than is typical for classical AMG or rAMG.  The
``Aggregation Operator" is then a binary matrix where column $i$ corresponds to
aggregate $i$ (i.e., the $i$th C-point) and this column is nonzero only for the
degrees-of-freedom in aggregate $i$.  See also \cite{root-node} which constructs initial sparsity patterns in this manner.} is used for $T$. This base pattern $T$ is then expanded 
based on the strength of connection matrix $S$ via $m-1$ multiplications in line \ref{alg:expandsparsity}.
We note that basing the interpolation pattern on strong connections is the same 
strategy as used by $\ell$AIR and root-node, and this approach allows us to generate nearly identical patterns.

The least squares solution for line \ref{alg:constraint1} is computed using a psuedoinverse based on $B_c$ restricted to the sparsity pattern of row $i$, $w^{\{i\},t}$. 
These pseudoinverses can be locally precomputed for efficiency and are typically small.

The projection operation in line \ref{alg:constraint2} takes each new update
$\witerupdate{*}{k}$ and projects it so that $Q \witerupdate{*}{k} = 0$.  That is, this
operation ensures that each update $\witerupdate{*}{k}$ does not disturb the mode
interpolation relationship $Q \witer{*}{t} = \B{*}|_F$.  The projection operation
with $Q$ can be implemented locally using the same strategy as for line \ref{alg:constraint1} (see \cite{OlScTu2011}).

\begin{algorithm}
\caption{\clair\ Algorithm}\label{alg:iterativeAIR}
\begin{algorithmic}[1]
   \STATE \textbf{Input}:\,$A$: Matrix
   \STATE $\phantom{\mbox{\textbf{Input}:}}$ $S$:  Strength matrix for interpolation
   \STATE $\phantom{\mbox{\textbf{Input}:}}$  $B$: User supplied mode constraint vector(s)
   \STATE $\phantom{\mbox{\textbf{Input}:}}$ $m$: Sparsity pattern degree 
   \STATE $\phantom{\mbox{\textbf{Input}:}}$ FC or Agg: Coarsening type
   \STATE \textbf{Output}: $P$: Interpolation in the form of the weight block \witer{*}{t} \vspace{4pt}
   \STATE \textbf{set} tentative prolongation \witer{*}{t}, corresponding to $[-A_{fc}, I]$
   \STATE \textbf{set} base sparsity pattern $T$ based on coarsening type
     \STATE $\phantom{**}$ \textbf{if} Coarsen type is FC
     \STATE $\phantom{*****}$ $T \leftarrow [S_{fc}, I]$
     \STATE $\phantom{**}$ \textbf{else if} Coarsen type is Agg
     \STATE $\phantom{*****}$ $T \leftarrow \mbox{Aggregation Operator}$\vspace{4pt}

   \STATE \textbf{set} expanded sparsity pattern to match F-row structure of $S^{m-1} T$ \label{alg:expandsparsity}\newline
    $\mathcal{W} \leftarrow \mbox{sparsity\_pattern}\left( (S^{m-1} T)|_F \right)$
   \STATE \textbf{expand} \witer{*}{t} to store (possibly zero) entries for every nonzero in $\mathcal{W}$ 
   \STATE \textbf{enforce} constraints on \witer{*}{t} such that  $Q \witer{*}{t} = \B{*}|_F$, by taking row $i$, $w^{\{i\},t}$, and computing  \label{alg:constraint1}\newline
       $w^{\{i\},t} \leftarrow$ least squares solution to $w^{\{i\},t} B_c = B$ \vspace{5pt}

   \STATE \textbf{compute} exact or approximate block inverses for $\widehat{\affinv{*}}$ \label{alg:invcompute}
   \FOR{$k = 1, 2, ...$}
   \STATE 
          $\phantom{*}$ \textbf{compute residual update} to minimize equation (\ref{eqn:minproba}), using $\mathbf{0}$ initial guess for $\witerupdate{*}{k}$\,:
          $$\witerupdate{*}{k} \leftarrow \widehat{\affinv{*}} \left( \left(-\afc{*} - \aff{*}\, \witer{*}{t} \right) - \aff{*} \mathbf{0} \right) $$
   \STATE $\phantom{*}$ \textbf{project update}: $\witerupdate{*}{k} \leftarrow (I - Q^T (Q Q^T)^{-1} Q) \witerupdate{*}{k}$ \label{alg:constraint2}
   \STATE $\phantom{*}$ \textbf{update interp}: $\witer{*}{t} \leftarrow \witer{*}{t} + \witerupdate{*}{k}$ \label{alg:finalupdate}
   \ENDFOR
   \STATE \textbf{set} final interpolation weights: $\w{*} \leftarrow \witer{*}{t}$
   \STATE \textbf{Return} \w{*}
\end{algorithmic}
\end{algorithm}

\begin{remark}
    If the constraints are removed from Algorithm \ref{alg:iterativeAIR} and the exact inverse $\affinv{*}$ is used, then $\ell$AIR is recovered.  That is, removing the constraint lines \ref{alg:constraint1} and
    \ref{alg:constraint2} and assuming $k=1$ yields  the final update in line \ref{alg:finalupdate} of
    \begin{align*}
        \witer{*}{t} &\leftarrow \witer{*}{t} + \witer{*}{1}\\
                     & \leftarrow \witer{*}{t} + \affinv{*} \left(-\afc{*} - \aff{*}\, \witer{*}{t} \right) \\
                     & \leftarrow - \affinv{*} \afc{*}
    \end{align*}
    However, we always use the constraints in Algorithm \ref{alg:iterativeAIR}, making \clair\ distinct from $\ell$AIR.
\end{remark}

\begin{remark}
   We will most often use the exact inverse $\affinv{*}$ for
   $\widehat{\affinv{*}}$ in line \ref{alg:invcompute}, similar to $\ell$AIR.
   In this case, Algorithm \ref{alg:iterativeAIR} is run with $k=1$, i.e.,
   the output of the algorithm does not change for $k>1$.
\end{remark}

\subsubsection{Comparison to Root-node}

We now clearly distinguish the similarities and differences between \clair\ and root-node. The \clair\ approach shares with root-node (i) the same mode and sparsity constraints (\ref{eqn:global_const1}) and (\ref{eqn:sparsityconst}), (ii) the ability to iteratively find $P$, and (iii) approximates some form of $\pideal$.

Regarding differences, root-node uses a simple diagonal approximation to $A_{ff}^{-1}$, whereas \clair\ uses a potentially much more powerful approximation $\widehat{\affinv{*}}$, where each block is often inverted exactly\footnote{Note that \clair\ also supports using a diagonal inverse similar to root-node, where each block of $\widehat{\affinv{*}}$ would be the local diagonal inverse.  However, this approximation does not always lead to effective AMG hierarchies in our experiments, e.g., for the 3D Poisson problem. Developing more effective approximations is future research.}. The proposed \clair\ method also uses a simple one-step iteration that is cheaper computationally than root-node, with no required global storage of Krylov vectors or the additional computations and communications required to maintain the (globally) orthogonal Krylov basis. Lastly, root-node has only considered fast aggregation-based coarsening, whereas \clair\ will support and explore both fast and slow coarsening, targeting diffusive and advective problems, respectively.

\section{Numerical Results}    
\label{sec:num}

We now present supporting numerical results for \clair, 
comparing against root-node and $\ell$AIR for a variety of classic 
model diffusion problems and for advection-diffusion problems over a range of diffusion parameters.  The proposed solver is shown to (i) achieve more desirable operator complexities than $\ell$AIR, (ii) compare well against the root-node solver for symmetric diffusion problems, where root-node is already known to perform well, and (iii) be more robust regarding parameter tuning when compared to $\ell$AIR, e.g., for the advection-diffusion problems when going from the purely advective to highly diffusive regimes.

Our basic test framework is as follows.  The $\ell$AIR and root-node solvers use the library implementations in PyAMG \cite{BeOlScSo2023} and \clair\ is also implemented in PyAMG\footnote{The \clair\ implementation is in the \texttt{CF\textunderscore rootnode} branch, commit bf56195b55a27bd99625c385ee9ba1df8814f764.}.  We use V(1,1) cycles as a preconditioner for CG and GMRES in the symmetric and nonsymmetric settings, respectively. The smoother choices respect the reduction framework.  For presmoothing, we use 1 iteration of CFF-weighted-Jacobi, i.e., one Jacobi sweep on the C-points, followed by two Jacobi sweeps on the F-points.  The weight is equal to $1/\rho(D^{-1}A)$. (Here, $\rho(\cdot)$ denotes the spectral radius, which is approximated numerically by 15 iterations of Arnoldi.) For postsmoothing, we use 1 iteration of FFC-weighted-Jacobi, which is defined analogously to CFF. Such relaxation methods are chosen for simplicity, parallelism, and the preservation of a symmetric preconditioner when $A$ is symmetric, so that CG can be used.
The absolute halting criteria is $10^{-9}$ for the smallest problem size, and is then scaled to simulate a discrete L2-norm\footnote{In 2D, the tolerance is scaled by 2 each grid refinement, and in 3D, by $\sqrt{8}$.}. We report preconditioned Krylov iterations and work-per-digit of accuracy.  One work unit is defined to be the cost of a fine-grid matrix-vector product, and work-per-digit of accuracy estimates how many work units are required to reduce the residual by one order of magnitude.  Work measurements allow for easy comparison across methods.  To derive our work-per-digit measure, we first estimate the total cost of doing one matrix-vector product at each level in the hierarchy, relative to one matrix-vector product on the finest level.  This yields the operator complexity measure
\begin{equation}
\mbox{OC} = \sum_k \mbox{nnz}(A_k) / \mbox{nnz}(A_0),
\end{equation}
where $A_k$ is the operator on level $k$ in the multigrid hierarchy and  $\mbox{nnz}(\cdot)$ stands for the number of nonzero entries.  To account for the cost of doing relaxation at each level, we multiply OC by 3.5 because we estimate (roughly) the cost of CFF- and FFC-Jacobi to be slightly less than 4 matrix-vector products. Thus work-per-digit of accuracy is estimated as
\begin{equation}
    \mbox{wpd} = 3.5\; \mbox{OC} / \log_{10}(\gamma),
\end{equation}
where the average residual convergence factor is $\gamma = (\| r_k\| / \|r_0\|)^{1/k}$, $r_k$ is the final residual, and $r_0$ is the initial residual.

\subsection{Diffusion Tests}
\label{sec:difftests}

For the symmetric case, we consider a variety of classic diffusion tests.
\subsubsection*{2D Poisson} The PDE here is $-u_{xx} - u_{yy} = f$, with Dirichlet boundary conditions on the unit box. The discretization is classic second-order 5-point finite differencing on a regular grid.

\subsubsection*{3D Poisson}
The PDE here is $-u_{xx} - u_{yy} - u_{zz} = f$, with Dirichlet boundary conditions on the unit box. The discretization is classic second-order 7-point finite differencing on a regular grid.

\subsubsection*{2D Grid-Aligned Anisotropic Diffusion} The PDE here is 
\begin{subequations}
\begin{align}
-\nabla \cdot Q^T D Q \nabla u &= f \quad \mbox{for } \Omega = [0,1]^2, \label{eqn:anisodiff1}\\
u &= 0 \quad \mbox{on } \partial \Omega, \label{eqn:anisodiff2}
\end{align}
\end{subequations}
where $\Omega$ is the unit box domain and
\begin{equation} 
Q = \begin{bmatrix} \cos(\varphi) & -\sin(\varphi)\\ \sin(\varphi) & \cos(\varphi)\end{bmatrix}, \quad D = \begin{bmatrix} 1 & 0\\ 0 & \epsilon \end{bmatrix}
\end{equation}
represent the rotation by angle $\varphi$ and the strength of anisotropy $\epsilon$.  The discretization uses a regular grid and bilinear (Q1) finite elements.  For this problem $\varphi = 0$ and $\epsilon = 0.001$, representing strong grid-aligned anisotropy in the $x$-direction.

\subsubsection*{2D Rotated Anisotropic Diffusion}
Here the PDE and discretization are the same as for equations (\ref{eqn:anisodiff1})--(\ref{eqn:anisodiff2}), but $\varphi = \pi/8$ and $\epsilon$ remains 0.001.  This represents strong non-grid-aligned anisotropy at the angle of $\pi/8$.  As noted in \cite{Jacob,KahlRott2018,BranBranKahlLivs2015b}, this is a difficult discretization and angle of anisotropy for AMG.

\subsubsection*{Box-in-Box Coefficient Jump} The PDE here is 
\begin{subequations}
\begin{align}
-\nabla \cdot d(x,y) \nabla u &= f \quad \mbox{for } \Omega, \label{eqn:jumpy1}\\
u &= 0 \quad \mbox{on } \partial \Omega, \label{eqn:jumpy2}
\end{align}
\end{subequations}
where $d(x,y)$ is the jumping coefficient.  Here, $\Omega = [0,1]^2$, $d(x,y) = 1 \mbox{ if } (x,y) \notin [0.44,0.52]^2$, and
$d(x,y) = 10^4 \mbox{ if } (x,y) \in [0.44,0.52]^2$.  The grid is regular and the coefficient jumps are grid-aligned on the finest level, but will not be grid-aligned at coarser levels due to the algebraic coarsening.  The discretization is classic second-order 5-point finite differencing from \cite{alcouffe1981multi} for coefficient jump problems. 

\subsubsection*{Sawtooth Coefficient Jump}
\begin{figure}[h!]
     \centering
         \includegraphics[width=0.32\textwidth]{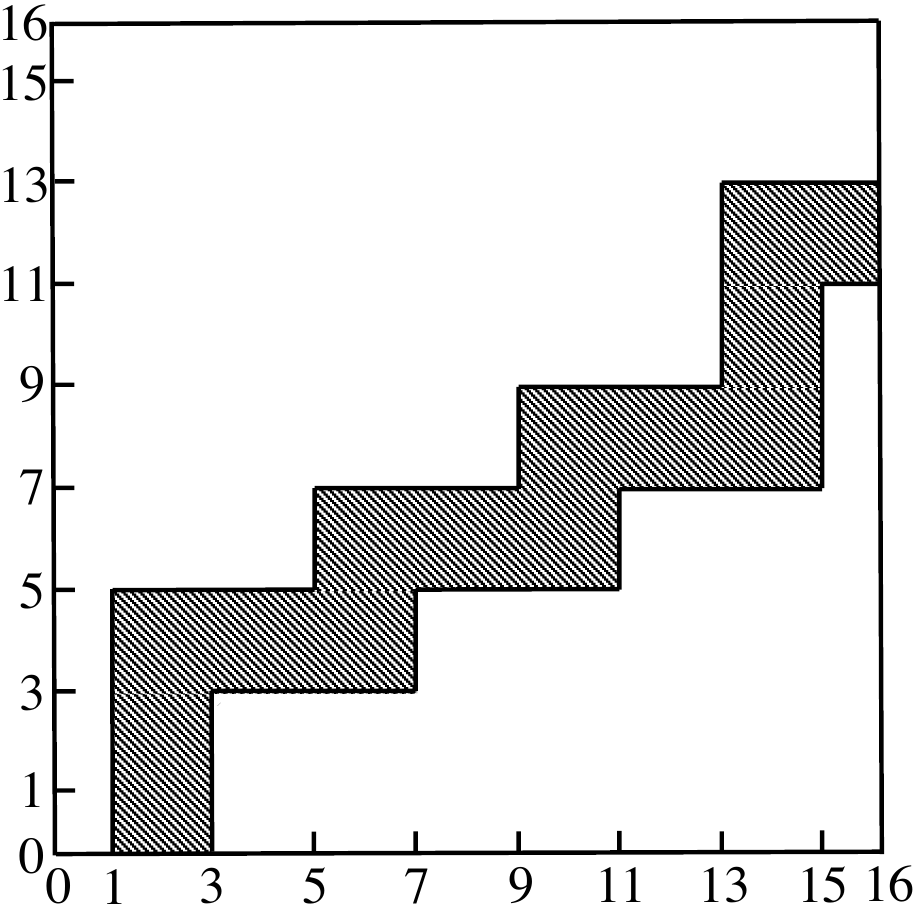}
         \caption{Sawtooth coefficient jump domain, $d(x,y)=10^4$ in shaded region and $d(x,y)=1$ outside.}
         \label{fig:sawtooth}
\end{figure}

Here, the PDE and discretization are the same as for equations (\ref{eqn:jumpy1})--(\ref{eqn:jumpy2}), but $\Omega = [0,16]^2$,  $d(x,y) = 1$ for points outside the shaded region of 
Figure \ref{fig:sawtooth}, and $d(x,y) = 10^4$ for points inside the shaded region \cite{alcouffe1981multi}.
\subsubsection*{Laplace Problem with Adaptive Mesh Refinement (AMR)} 
\begin{figure}[h!]
     \centering
     \begin{subfigure}[b]{0.29\textwidth}
         \centering         \includegraphics[width=\textwidth]{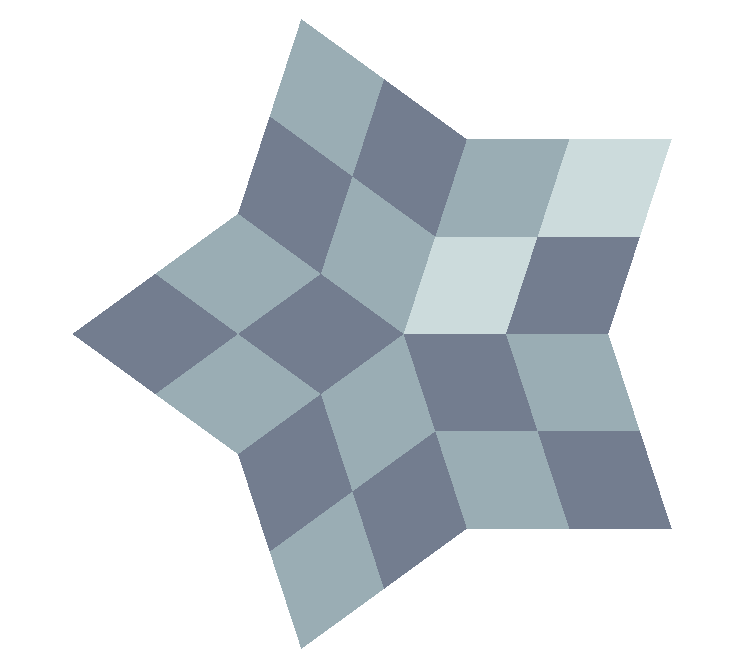}
         \caption{Initial star mesh }
         \label{fig:mesh_1}
     \end{subfigure}
     \hfill
     \begin{subfigure}[b]{0.29\textwidth}
         \centering
         \includegraphics[width=\textwidth]{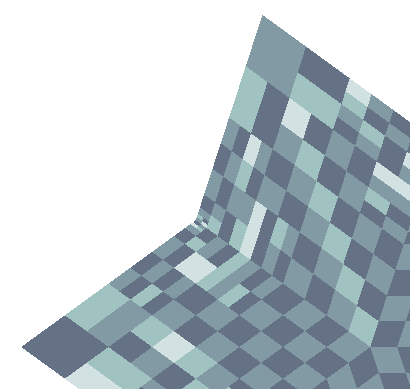}
         \caption{Re-entrant corner after $5^{\text{th}}$ refinement}
         \label{fig:mesh_2}
     \end{subfigure}
     \hfill
     \begin{subfigure}[b]{0.29\textwidth}
         \centering
         \includegraphics[width=\textwidth]{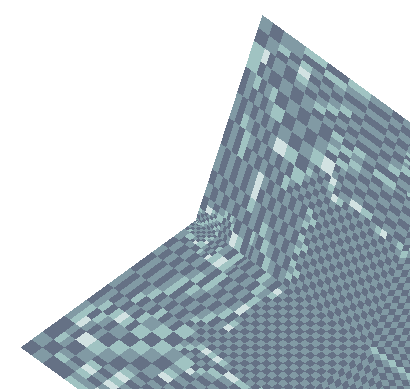}
         \caption{Re-entrant corner after $10^{\text{th}}$ refinement}
         \label{fig:mesh_3}
     \end{subfigure}
        \caption{AMR throughout different stages for the Laplace problem.}
        \label{fig:mesh_plots}
\end{figure}

Here, we consider a finite element discretization of the Laplace problem $-\Delta u =1$ with homogeneous Dirichlet boundary conditions. This problem  is solved on a sequence of meshes in an $\mathbb{H}^{1}$-conforming finite element space which are locally refined in accordance with a simple Zienkiewicz-Zhu \cite{zz_error} error estimator. We save the stiffness matrix $A$ throughout each iteration of the AMR loop and use this matrix to further illustrate differences in outcomes across various AMG solvers. This example is identical to PyMFEM library example 6 \cite{mfem,pymfem}. AMR stages for the star mesh file used in this problem are shown in Figure \ref{fig:mesh_plots}.
\subsubsection*{Standard Solver Parameters}
We now list our standard solver parameters for all three methods, $\ell$AIR, \clair, and root-node.  We will occasionally tune the parameters (e.g., strength) for $\ell$AIR and root-node, in order to make the existing methods more competitive.  The parameters for the proposed \clair\ method will remain fixed for all symmetric test cases, except for the 2D Rotated Anisotropic Diffusion case, where we will follow the guidance of \cite{Jacob} and consider larger degree interpolation sparsity patterns. We use accelerated CG for \clair\ and root-node and accelerated GMRES for $\ell$AIR, as the typical approach for $\ell$AIR is not symmetric (i.e., does not satisfy $R = P^T$, see below parameter choices).

\vspace{6pt} \noindent \textit{Parameters for $\ell$AIR} (see \cite{Ben} for more details on parameter definitions)
\squishlist
    \item Strength-of-connection tolerance for coarsening 0.25 
    \item Degree 2 $\ell$AIR restriction (sparsity pattern similar to $m=2$ in Algorithm \ref{alg:iterativeAIR}) with  interpolation strength tolerance 0.05
    \item Ruge-St\"{u}ben coarsening first-pass only~\cite{JWRuge_KStuben_1987a} (coarsen type FC in Algorithm \ref{alg:iterativeAIR}, relatively slow coarsening often used by $\ell$AIR)
    \item Coarse-grid matrices filtered with drop-tolerance of $10^{-4}$ 
    \item Classical interpolation formula used for $P$ \cite{JWRuge_KStuben_1987a}, which is typical for $\ell$AIR and diffusive problems
\squishend

\vspace{6pt} \noindent \textit{Parameters for \clair}
\squishlist
    \item Strength-of-connection tolerance for coarsening 0.5 
    \item Degree 2 sparsity pattern for $P$ ($m=2$ in Algorithm \ref{alg:iterativeAIR}) with interpolation strength tolerance 0.5 used to generate $S$ in Algorithm \ref{alg:iterativeAIR}
    \item Mode constraint vector $B = \mathbf{1}$, presmoothed with 5 iterations of CFF-weighted-Jacobi at each level
    \item $R = P^T$ used (see discussion below)
    \item Aggregation based coarsening (coarsen type Agg in Algorithm \ref{alg:iterativeAIR})
\squishend
 
\vspace{6pt} \noindent \textit{Parameters for root-node} (see \cite{root-node} for more parameter details)
\squishlist
    \item Strength-of-connection tolerance for coarsening 0.25
    \item $P$ smoothed with energy-minimizing projected CG and a degree 2 sparsity pattern (sparsity pattern similar to $m=2$ in Algorithm \ref{alg:iterativeAIR})
    \item Mode constraint vector $B = \mathbf{1}$, presmoothed with 5 iterations of Jacobi at each level
    \item $R = P^T$ used (see discussion below)
    \item Aggregation based coarsening (standard coarsening for root-node) \vspace{6pt}
\squishend

We note that the interpolation strength-of-connection tolerance for \clair\ is different than that for $\ell$AIR (0.5 versus 0.05).  Both of these tolerance were tuned individually for each method.  We also note that \clair\ and root-node both use $R=P^T$, while $\ell$AIR uses classical interpolation. The use of classical interpolation (instead of $R^T$, the transpose of $\ell$AIR restriction) with $\ell$AIR is typical for diffusive problems, and additionally, the use of $R^T$ as interpolation leads to undesirably large operator complexities.

\subsubsection{Poisson Results}

\begin{figure}[h!]
     \centering
     \begin{subfigure}[b]{0.49\textwidth}
         \centering
         \includegraphics[width=\textwidth]{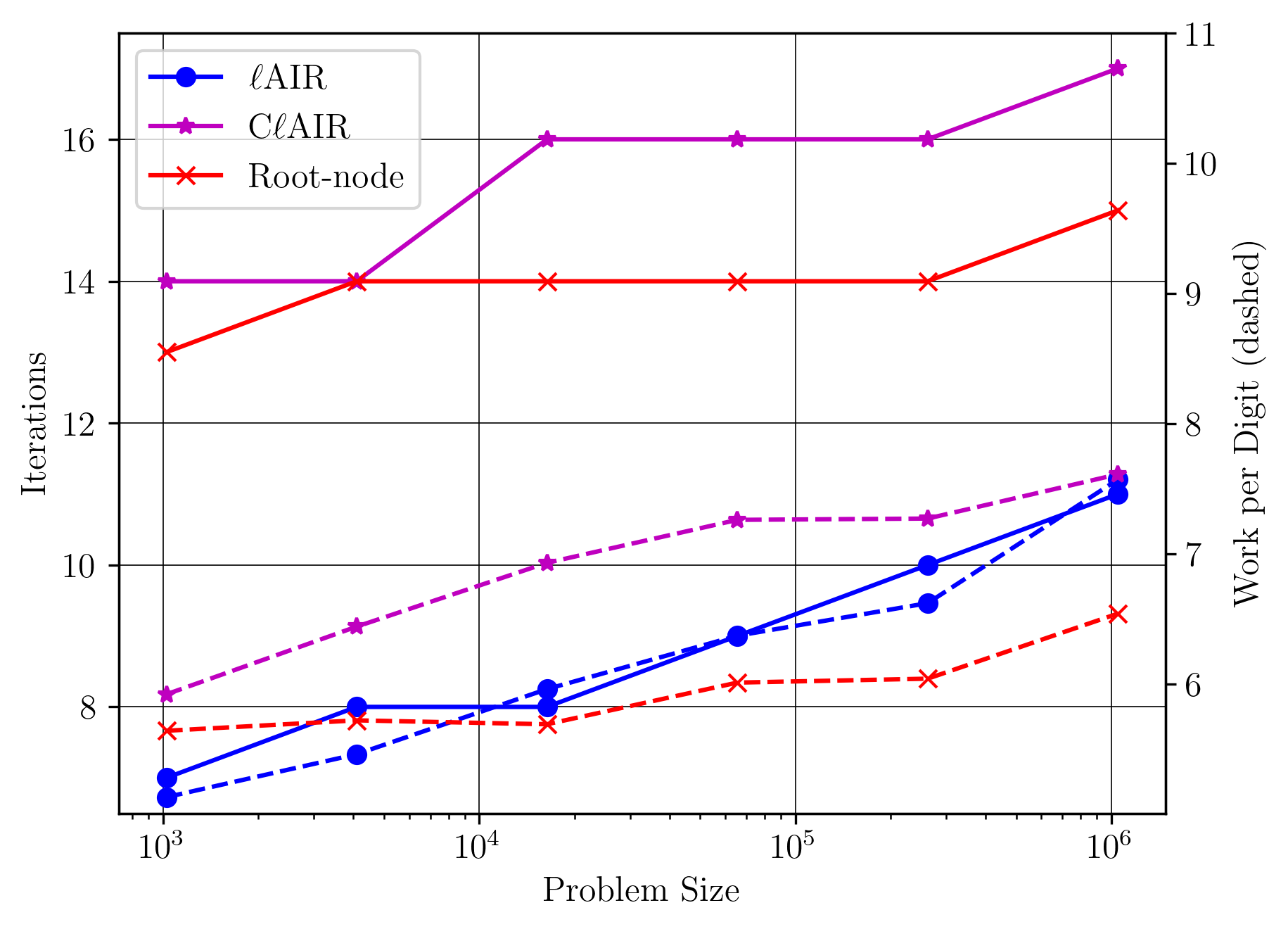}
         \caption{Poisson 2D }
         \label{fig:2d_poisson_1}
     \end{subfigure}
     \hfill
     \begin{subfigure}[b]{0.49\textwidth}
         \centering
         \includegraphics[width=\textwidth]{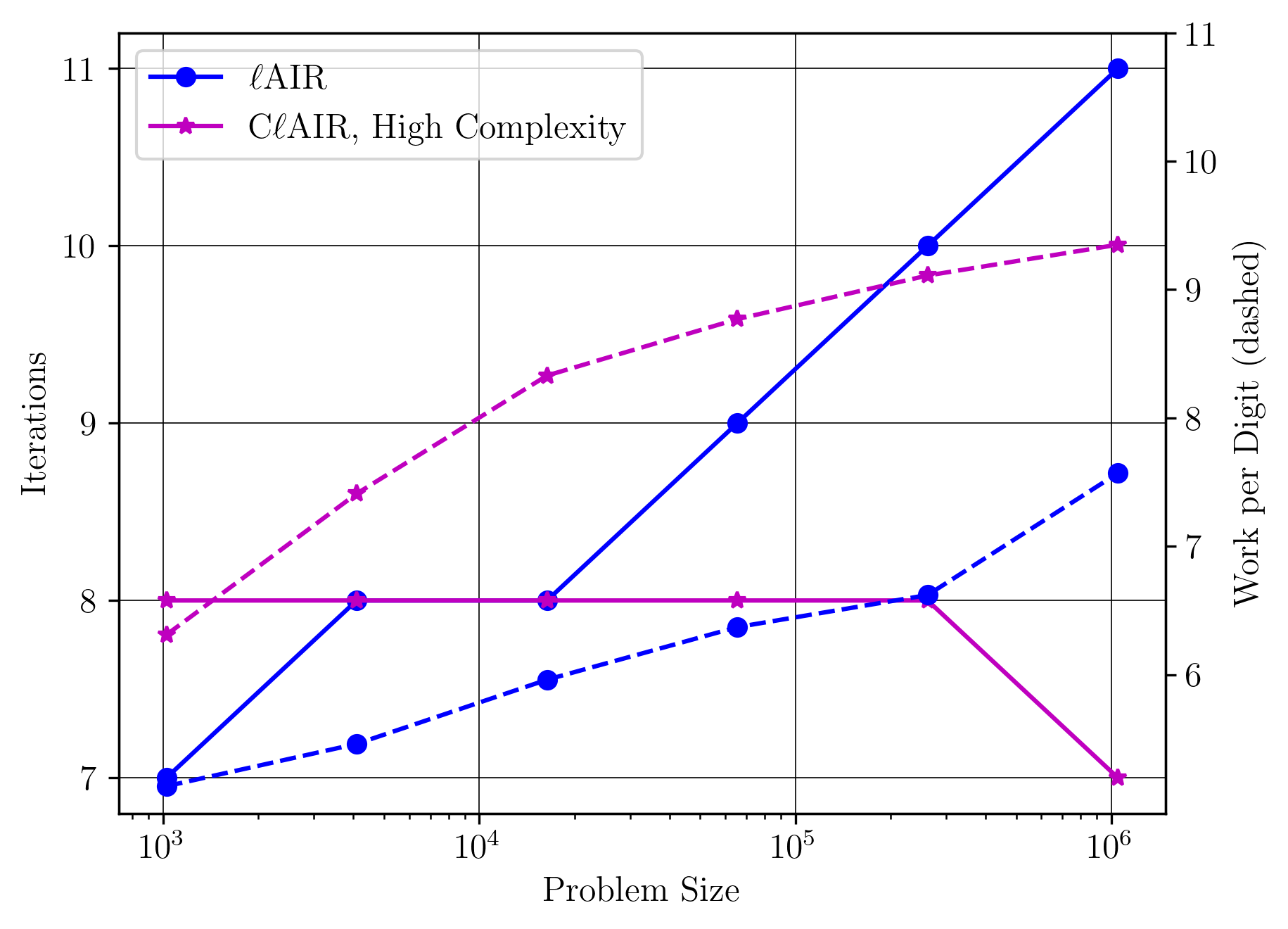}
         \caption{Poisson 2D, showing high complexity \clair}
         \label{fig:2d_poisson_2}
     \end{subfigure}
     \hfill
     \begin{subfigure}[b]{0.49\textwidth}
         \centering
         \includegraphics[width=\textwidth]{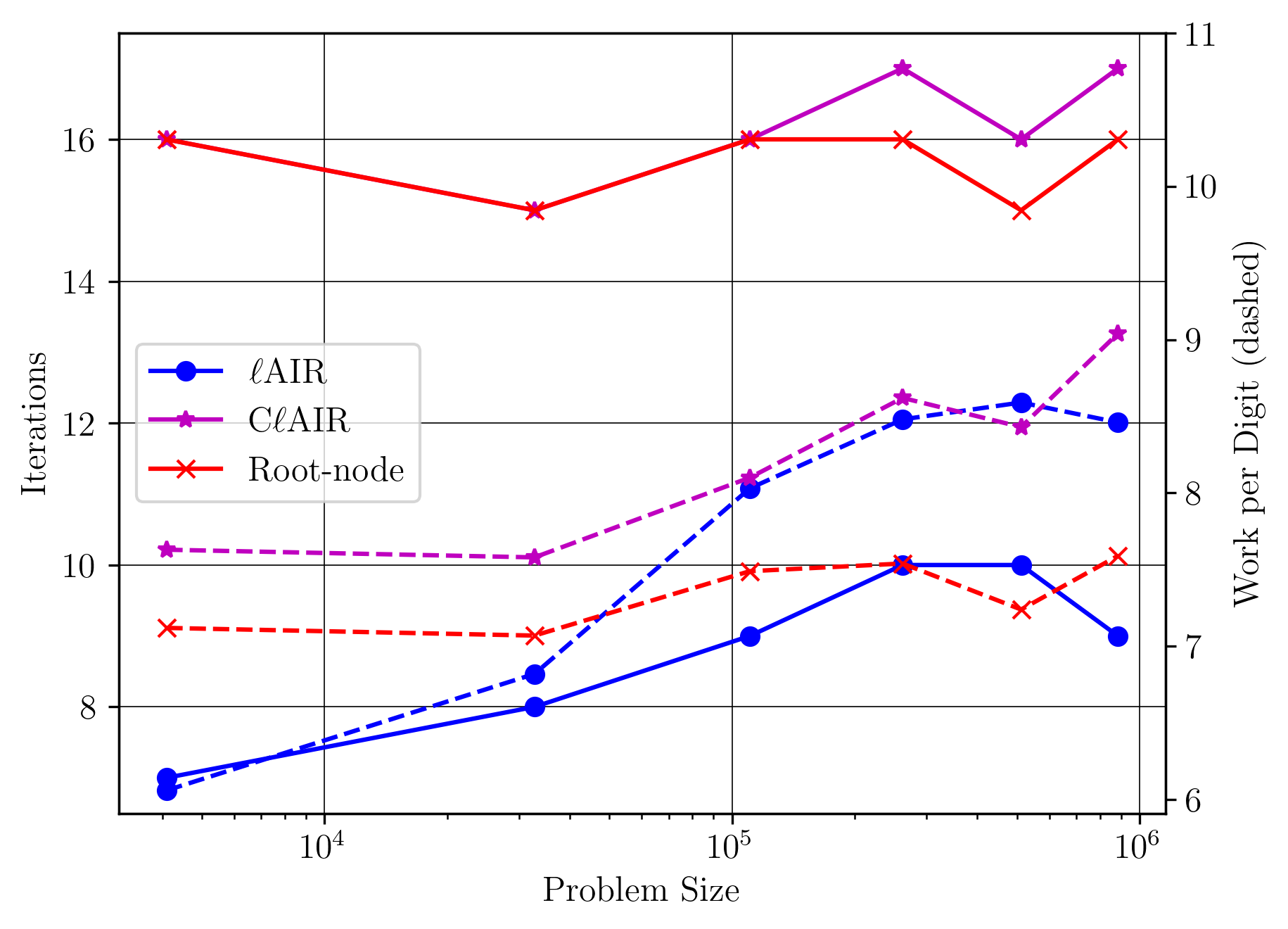}
         \caption{Poisson 3D}
         \label{fig:3d_poisson}
     \end{subfigure}
        \caption{2D and 3D Poisson, comparison of iterations and work-per-digit of accuracy for $\ell$AIR, \clair, and root-node.}
        \label{fig:Poisson_data}
\end{figure}

We next examine convergence results for the 2D and 3D Poisson problems in Figure \ref{fig:Poisson_data}.  Here, we tune $\ell$AIR to obtain a more challenging baseline solver and set the sparsity pattern in $R$ to degree 1.  For these simplest of problems, this choice reduces the operator complexity, while not affecting convergence.  For our other test problems, such a parameter choice does not lead to the best performance.  Operator complexities over all test problems, which are an advantage of \clair, are discussed later in Section \ref{sec:complexity}.

Figures \ref{fig:2d_poisson_1} and \ref{fig:3d_poisson} depict results for the 2D and 3D Poisson problems, respectively.  Overall, we see flat iteration counts for all methods, except a slight growth for $\ell$AIR in the 2D Poisson case. The chief performance difference is that $\ell$AIR has a significantly higher operator complexity than \clair\ and root-node, which will be discussed in Section \ref{sec:complexity}

Figure \ref{fig:2d_poisson_2} demonstrates that \clair\ can also be run similar to $\ell$AIR with Ruge-St\"{u}ben coarsening with first-pass only and classical interpolation for $P$.  In general, this setup does not lead to the most efficient solver for \clair, due to the high operator complexity ($\mbox{OC} \approx  3.3$).

\subsubsection{Other Diffusion Results}

\begin{figure}[h!]
     \centering
     \begin{subfigure}[b]{0.49\textwidth}
         \centering
         \includegraphics[width=\textwidth]{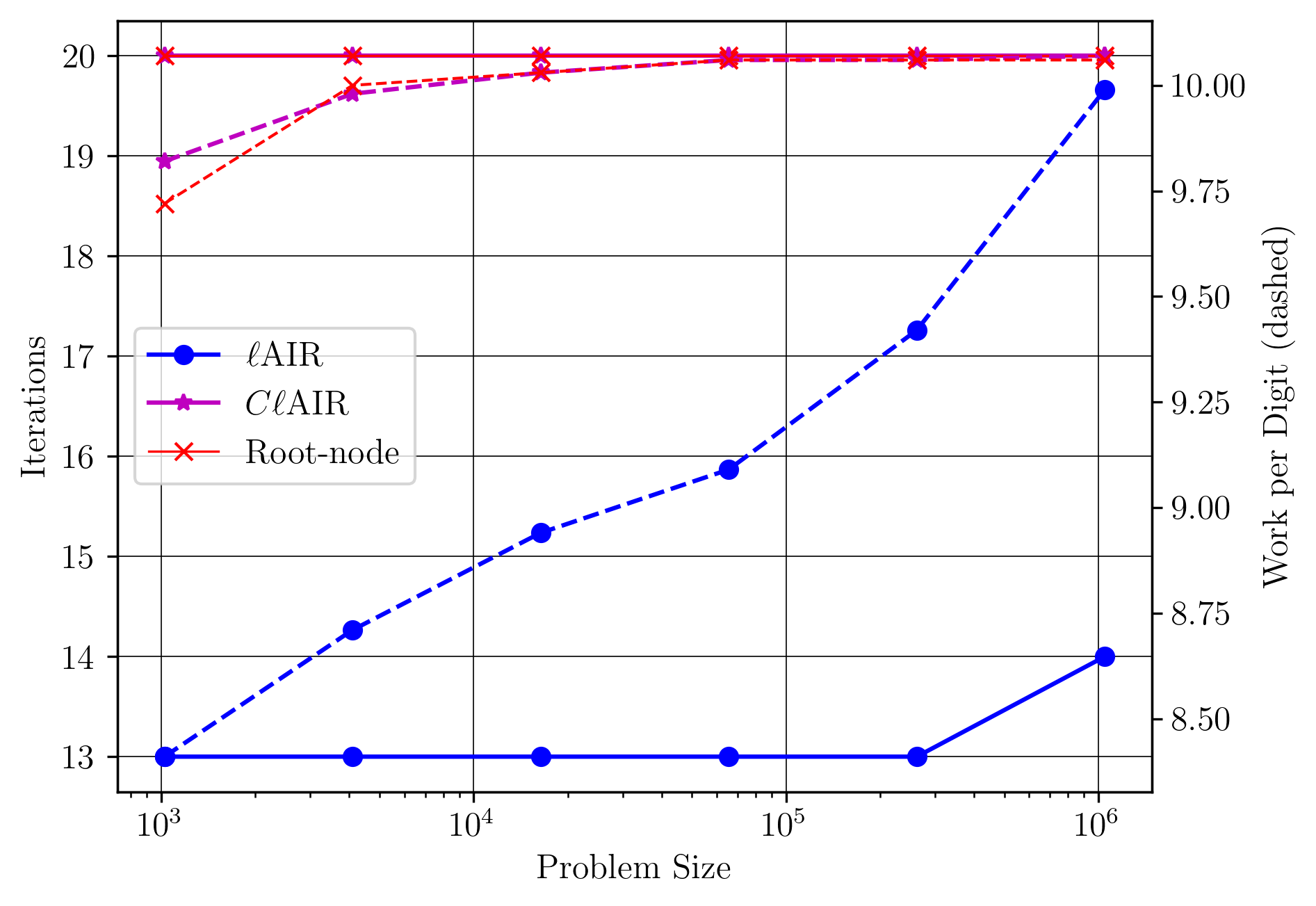}
         \caption{2D Grid-aligned Anisotropic Diffusion, when not visible, the lines for \clair\ are underneath root-node.}
         \label{fig:diffusion_1}
     \end{subfigure}
     \hfill
     \begin{subfigure}[b]{0.49\textwidth}
         \centering
         \includegraphics[width=\textwidth]{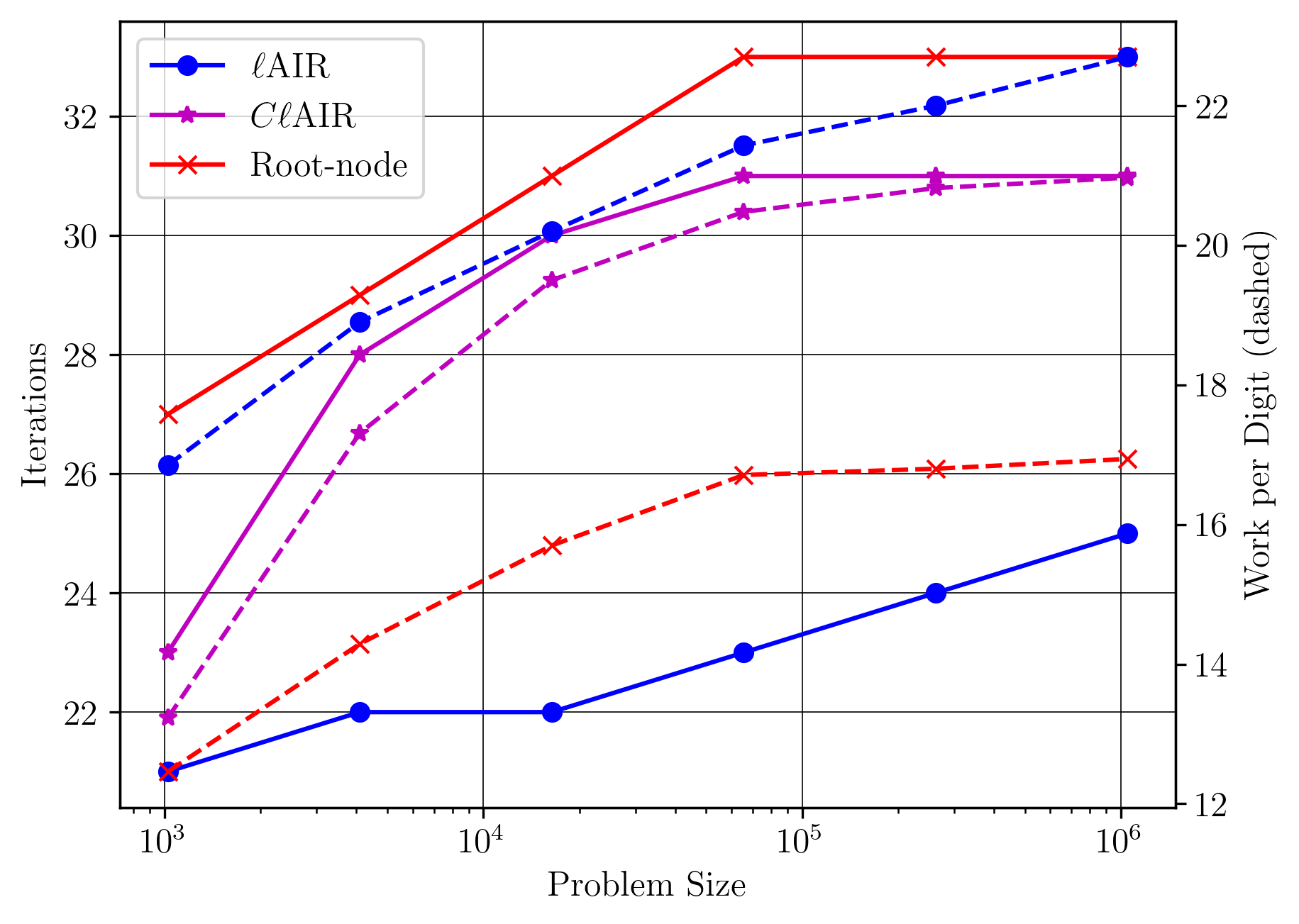}
         \caption{2D Rotated Anisotropic Diffusion  }
         \label{fig:diffusion_2}
     \end{subfigure}
     \hfill
     \begin{subfigure}[b]{0.49\textwidth}
         \centering
         \includegraphics[width=\textwidth]{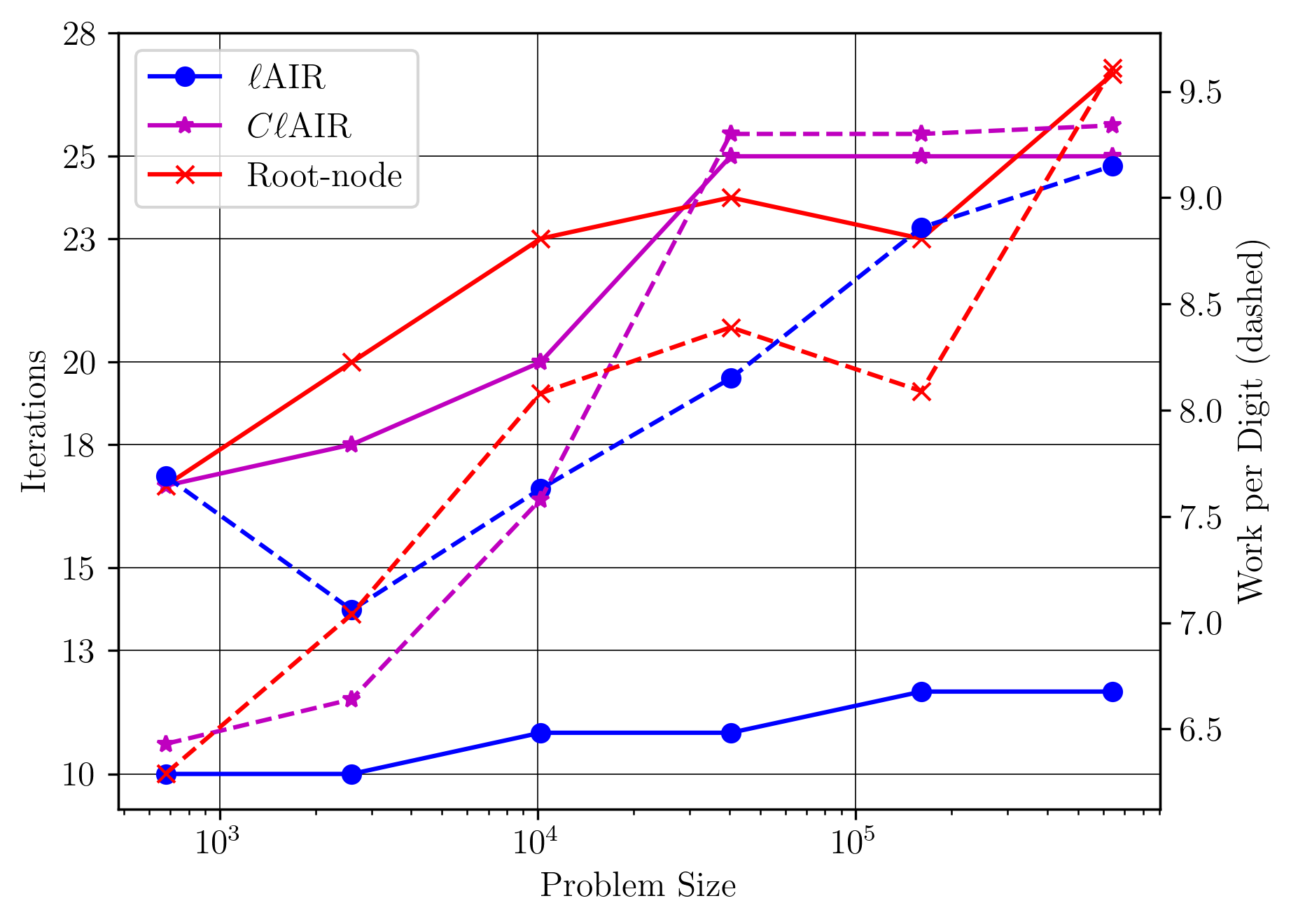}
         \caption{Box-in-Box Coefficient Jump}
         \label{fig:diffusion_3}
     \end{subfigure}
     \hfill
     \begin{subfigure}[b]{0.49\textwidth}
         \centering
         \includegraphics[width=\textwidth]{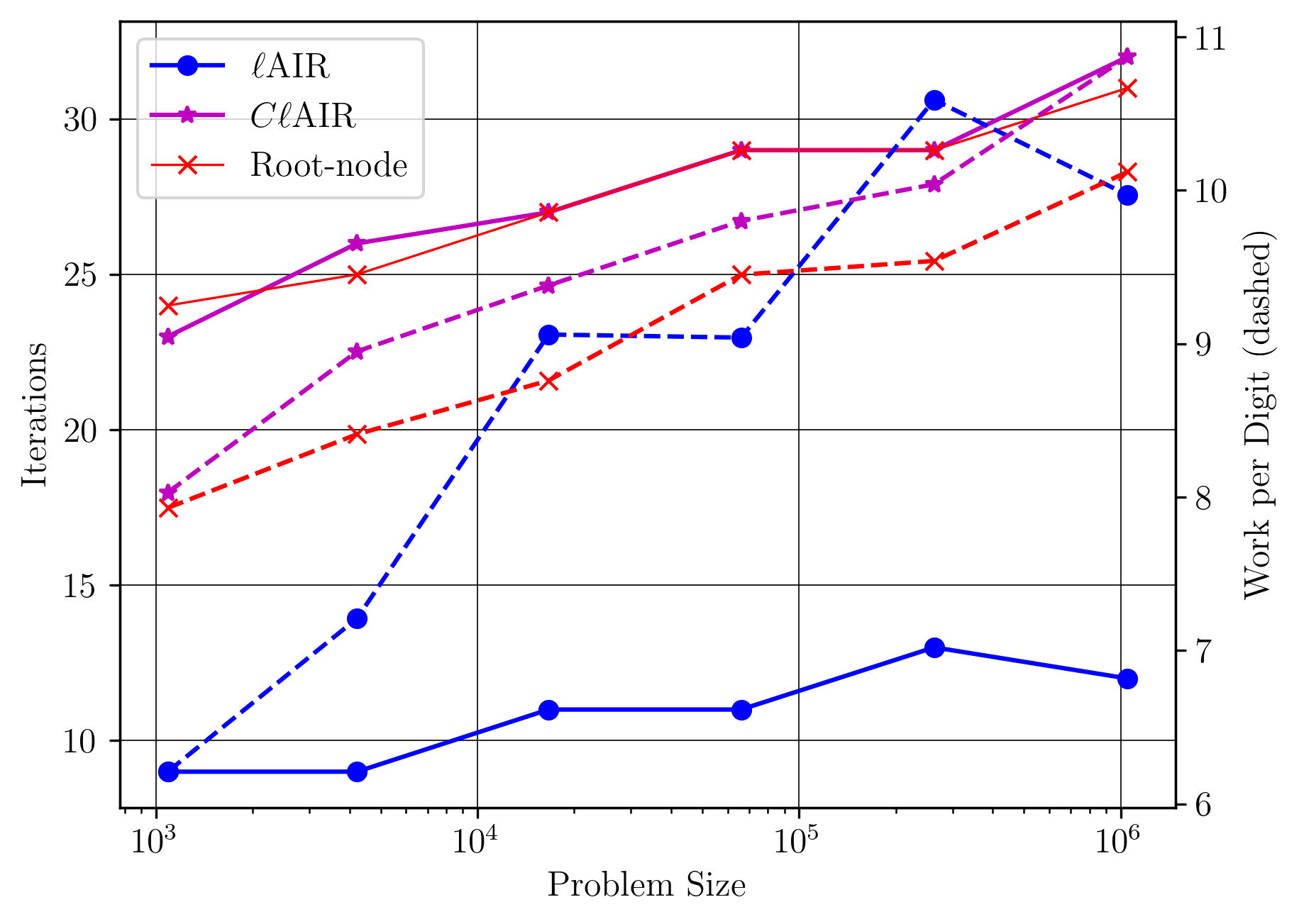}
         \caption{Sawtooth Coefficient Jump, when not visible, the lines for \clair\ are underneath root-node.}
         \label{fig:diffusion_4}
     \end{subfigure}
     \hfill
     \begin{subfigure}[b]{0.49\textwidth}
         \centering
         \includegraphics[width=\textwidth]{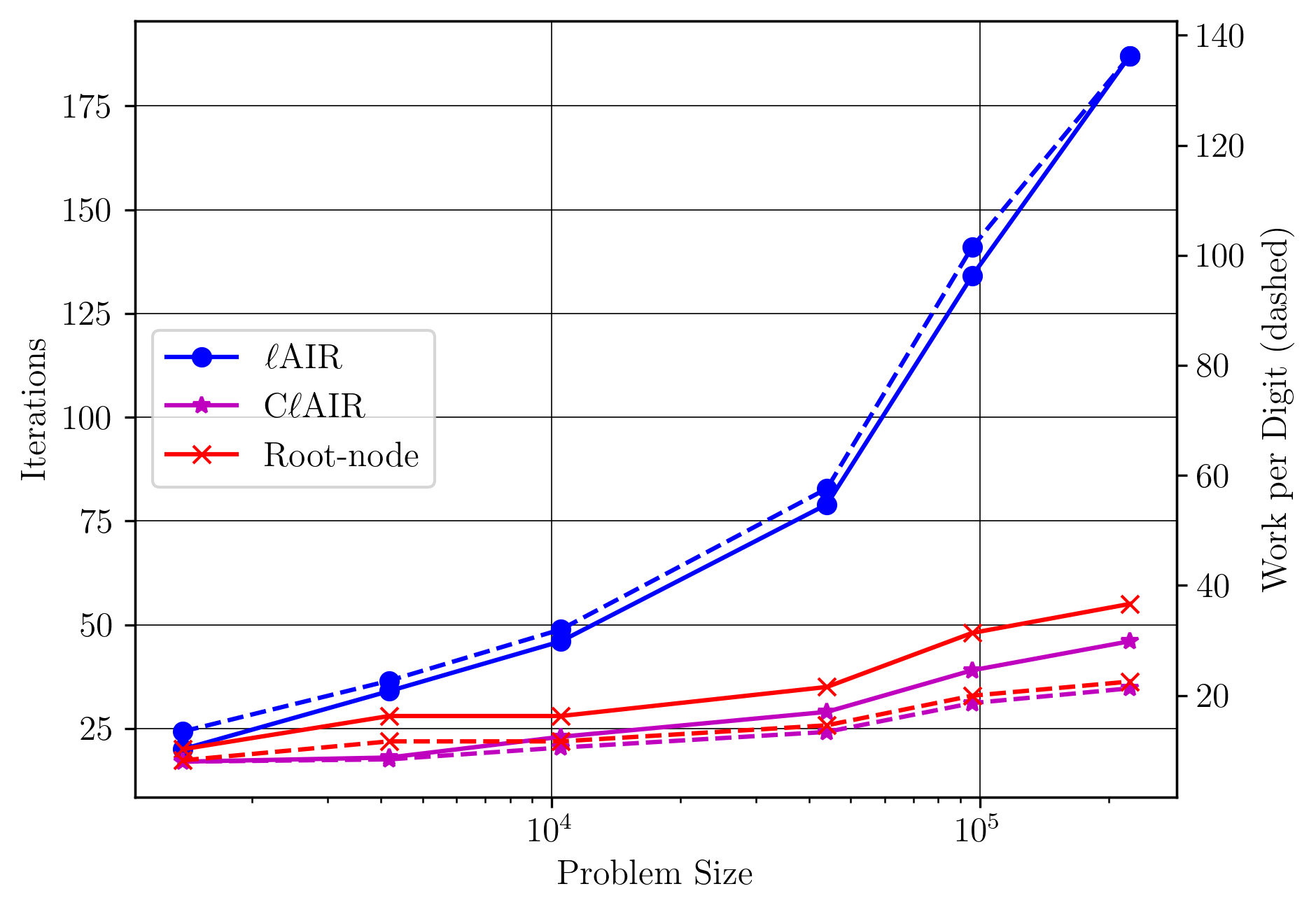}
         \caption{Laplace Problem with AMR.}
         \label{fig:diffusion_5}
     \end{subfigure}
        \caption{Various diffusion test problems, comparison of iterations and work-per-digit of accuracy for $\ell$AIR, \clair, and root-node.}
        \label{fig:diffusion_data}
\end{figure}

We now consider our other diffusion test problems in Figure \ref{fig:diffusion_data}, which again depicts iterations and work-per-digit of accuracy. For these tests, we tune root-node and $\ell$AIR on the 2D Grid-Aligned Anisotropic Diffusion problem, in order to obtain more challenging baseline solvers.  Here, we set the coarsening and interpolation strength-of-connection tolerances to 0.5 and note that using such a high interpolation strength tolerance typically hurts $\ell$AIR performance, with this case being the outlier.  Thus, this is not a general parameter setting for $\ell$AIR and is used to highlight the greater flexibility of the untuned \clair\ solver. Following the work \cite{Jacob}, we also consider higher-degree interpolation for the 2D Rotated Anisotropic Diffusion case, for both root-node and \clair.  Unfortunately, prohibitive increases in operator complexity for $\ell$AIR do not allow for such an examination of higher-degree interpolation.

Figure \ref{fig:diffusion_1} depicts results for the 2D Grid-Aligned Anisotropic Diffusion problem, where we see similar performance for \clair, root-node, and $\ell$AIR. Figure \ref{fig:diffusion_2} depicts results for the 2D Rotated Anisotropic Diffusion problem where all three solvers produce similar work-per-digit numbers, and the root-node and \clair\ solvers show flat, scalable iteration counts for the higher-degree interpolation option.

Figures \ref{fig:diffusion_3} and \ref{fig:diffusion_4} depict results for the Box-in-Box Coefficient Jump and Sawtooth Coefficient Jump problems. All three solvers again show similar work-per-digit numbers, with the main difference being in operator complexity, which we discuss in the next subsection. 

Finally, results for the Laplace problem with AMR utilizing an unstructured star mesh are shown in Figure \ref{fig:diffusion_5}. While \clair\, and root-node performance are comparable in this instance, $\ell$AIR shows high preconditioned iterations and work-per-digit accuracy as problem size increases. $\ell$AIR performance may be improved by raising the strength-of-connection tolerance from $0.25$ to $0.5$ and adding the second pass of Ruge-St\"{u}ben coarsening, resulting in a maximum of $25$ preconditioned iterations. However, in doing so, the operator complexity becomes prohibitively high at around $3.88$. Lastly, we note that the work per digit for \clair\ and root-node grows slowly for this problem as the mesh is adaptively refined, but we do not find this surprising given the re-entrant corners and adaptive refinement.  We also observed a similar slow growth when using the benchmark  classical Ruge-St\"{u}ben solver for this problem.

\subsubsection{Operator Complexity Comparison}
\label{sec:complexity}

One key advantage for \clair\ when compared to other reduction-based approaches on diffusion problems is the ability to achieve moderate operator complexities with good convergence.  To illustrate this, Table \ref{table_op_comp} depicts the operator complexity for all three solvers and the largest problem size for each of the diffusion test problems.  Here, this advantage to \clair\ becomes obvious.  The faster aggregation-based coarsening, which only root-node and \clair\ support, allows for dramatically lower operator complexities and the associated lower storage requirements. 

Importantly, one could not simply introduce the aggregation-based coarsening into $\ell$AIR and maintain good convergence. The addition of the mode interpolation constraint vector $B$ is needed for good convergence.

We note that the operator complexity of 1.98 for \clair\ and the 2D Rotated Anisotropic Diffusion case can easily be brought down to 1.78 without a meaningful effect on convergence by using 0.25 as the strength-of-connection tolerance, but we choose to maintain uniform parameters for \clair.

\begin{table}[h!]
\begin{center}
\begin{tabular}{ c|c|c|c } 
\toprule
\textbf{Test Problem}                   & $\ell$AIR OC & \clair\ OC & Root-node OC \\
\midrule
2D Poisson                              & $2.20$ & $1.40$ & 1.34 \\ 
3D Poisson                              & $2.87$ & $1.71$ & 1.57 \\ 
2D Grid-aligned Anisotropic Diffusion   & $4.92$ & $1.52$ & 1.51\\
2D Rotated Anisotropic Diffusion        & $2.89$ & $1.98$ & 1.52 \\ 
Box-in-Box Coefficient Jump             & $2.73$ & $1.40$ & 1.34 \\ 
Sawtooth Coefficient Jump               & $2.71$ & $1.42$ & 1.35 \\
Laplace Problem with AMR              & $2.03$ & $1.32$ & 1.17 \\
\bottomrule
\end{tabular}
\end{center}
\caption{Operator complexities of $\ell$AIR and \clair\, for different diffusion test problems.}
\label{table_op_comp}
\end{table}

\begin{remark}
In summary, we have shown that \clair\ performs similarly in terms of work-per-digit, operator complexity, and iterations on this suite of classic diffusion test problems when compared to root-node, an AMG method known to be efficient for such problems.  Our goal is not to find a faster solver for the 5-point Poisson operator, which would be difficult and of questionable value, but instead to verify the proposed solver.  We have furthermore shown greater parameter insensitivity and dramatically lower operator complexities when compared to $\ell$AIR for these problems.  Finally, it is important to highlight that \clair\ retains good convergence and low operator complexity across a range of diffusion problems, which has proven challenging for previous reduction-based approaches. 
\end{remark}

\subsection{Advection-Diffusion Tests}
\label{sec:advtests}

We now turn our attention to nonsymmetric tests and consider the following two advection-diffusion  problems.

\subsubsection*{2D Constant Direction Advection-Diffusion}  Here the PDE is 
\begin{subequations}
\begin{align}
-\alpha \nabla \cdot \nabla u + \mathbf{b}(x,y) \cdot \nabla u &= f \quad \mbox{for } \Omega = [-1,1]^2, \label{eqn:advdiff1}\\
u &= 0 \quad \mbox{on } \partial \Omega, \label{eqn:advdiff2}
\end{align}
\end{subequations}
where $\alpha$ is the diffusion constant and $\mathbf{b}(x,y)$ describes the advection.  For this problem $\mathbf{b}(x,y) = [ \sqrt{2/3},\, \sqrt{1/3} ]$ representing constant non-grid-aligned advection.  The discretization is first-order upwinded discontinuous Galerkin (DG) for advection and the interior-penalty method for the DG discretization of diffusion.  The discretization was generating by utilizing examples 14 (diffusion) and 9 (advection) with PyMFEM \cite{mfem,pymfem}. For the case where $\alpha=0$, the boundary conditions are modified slightly such that an outflow now occurs on the North and East walls ($u=0$ is still prescribed on the West and South walls). We use $u=0$ for the inflow, so that we may test our solver with a zero right-hand-side and random initial guess, as is commonly done to verify AMG solvers.

\subsubsection*{2D Recirculating Advection-Diffusion} Here the PDE and discretization is the same as for equations (\ref{eqn:advdiff1})--(\ref{eqn:advdiff2}), 
except $\mathbf{b}(x,y)= [ x(1-x)(2y-1),\, -(2x-1)(1-y)y]$, representing divergence-free recirculating advection about the origin. For $\alpha=0$, this problem is ill-defined, thus we demonstrate numerical results over a range of diffusion constants ($\alpha$), starting from a smallest diffusion constant $\alpha=10^{-4}$ and going to $\alpha=10.0$, so that we test the advective, mixed, and diffusive regimes.

\subsubsection{Advection-Diffusion Results}

Our goal here is to test the solver's robustness from the highly advective to highly diffusive regime and show improved performance and robustness relative to the current state $\ell$AIR. Similar to \cite{Ben}, we pre-scale the fine-grid matrices with the inverse of their diagonal block (block-size equals the DG element size of 4).  Our solver parameters will remain fixed over these tests, and involve only minor changes to the symmetric parameters from Section \ref{sec:difftests}.  As the matrices are nonsymmetric, GMRES is accelerated with V(1,1)-cycles. The relaxation weight for postsmoothing with FFC-Jacobi is removed, as we no longer need symmetry for our preconditioner and this change slightly improves convergence for all methods.\footnote{The removal of the weight also more closely follows the reduction point-of-view in that post-smoothing should primarily solve the F-equations, so the slight improvement of convergence is not surprising.} The parameters for $\ell$AIR were changed slightly to use the second pass of Ruge-St\"{u}ben coarsening.  The parameters for \clair\ changed slightly with the strength-of-connection parameters becoming the same as for $\ell$AIR (0.25 for coarsening and 0.05 for interpolation).  Additionally because of the nonsymmetry, $P$ is generated separately from $R$ using $A$ and $A^T$, respectively. We also again consider two variants of \clair, a high complexity version that uses first pass only Ruge-St\"{u}ben coarsening, and lower complexity versions that uses aggregation-based coarsening. 
The parameters for root-node changed similarly, where GMRES is now used for the energy-minimization when computing $P$ and $P$ is generated separately from $R$ using $A$ and $A^T$, respectively.

\begin{figure}[h!]
     \centering
     \begin{subfigure}[b]{0.49\textwidth}
         \centering
         \includegraphics[width=\textwidth]{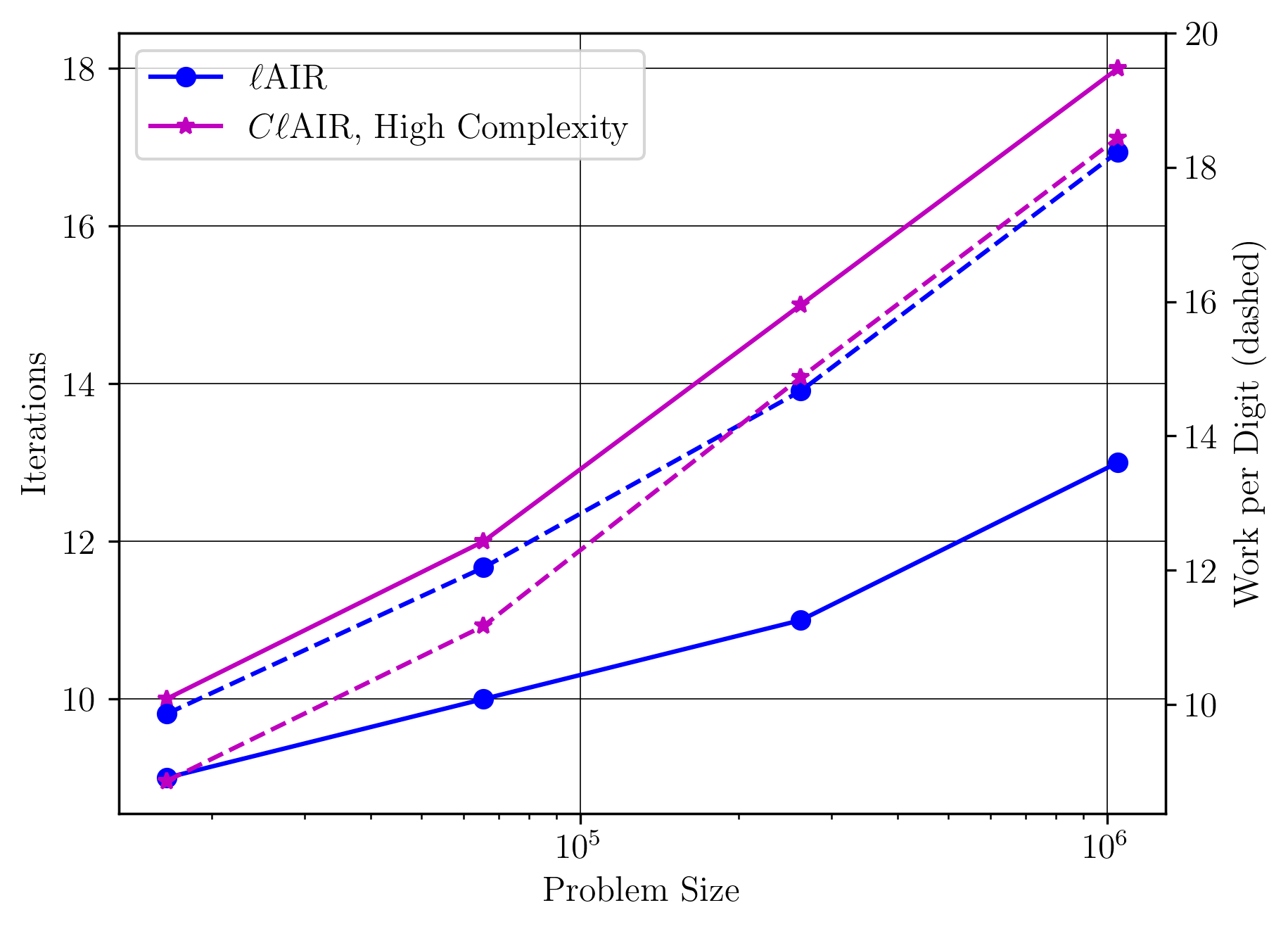}
         \caption{Solver performance for 2D Constant Direction Advection, but no diffusion ($\alpha=0$), comparing $\ell$AIR and \clair. Note, root-node does not converge within 100 iterations for this problem and is not shown.}
         \label{fig:adv_diff_1}
     \end{subfigure}
     \hfill
     \begin{subfigure}[b]{0.49\textwidth}
         \centering
         \includegraphics[width=\textwidth]{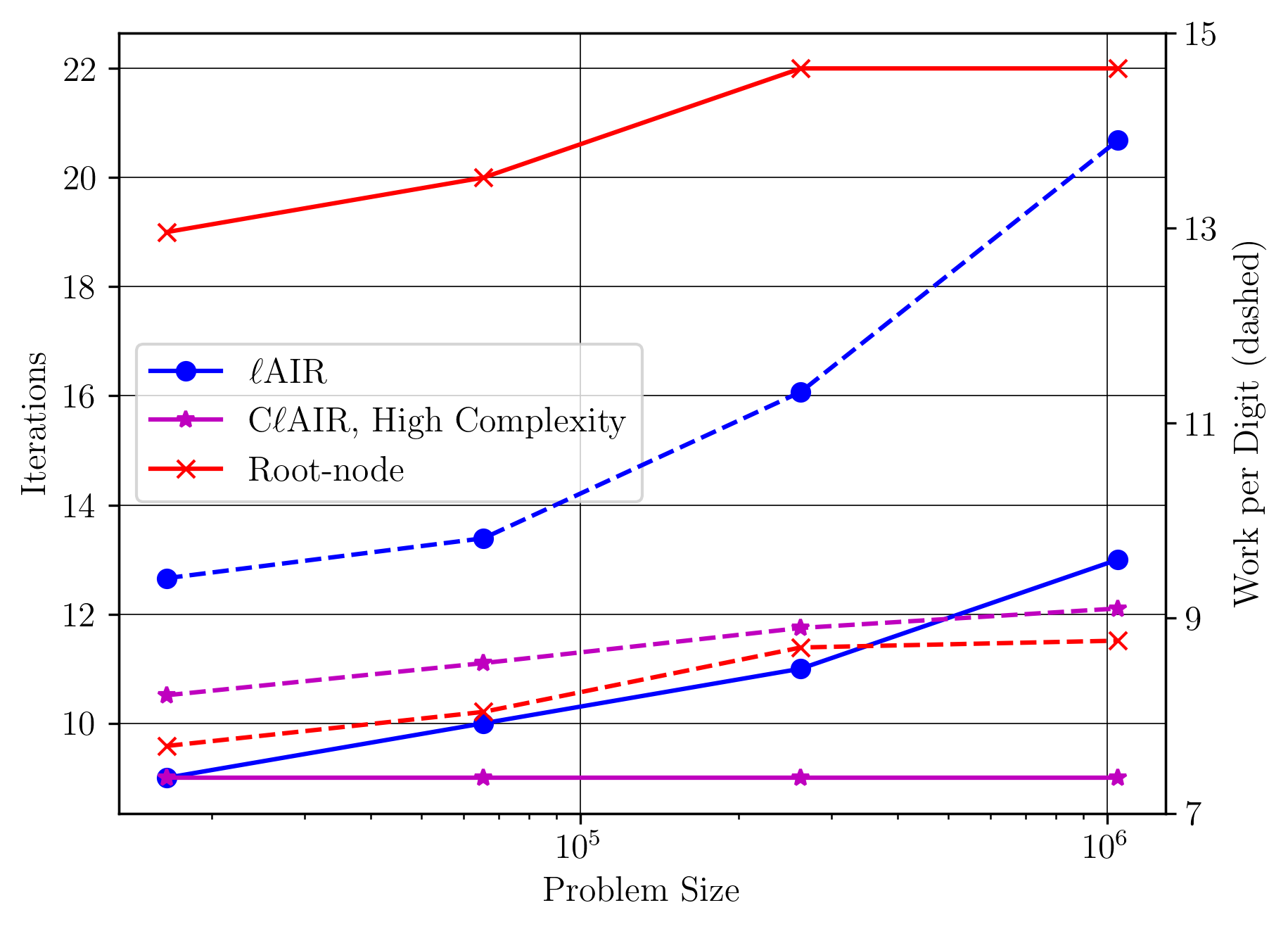}
         \caption{Solver performance for 2D Constant Direction Advection with diffusion ($\alpha=10$), comparing $\ell$AIR, \clair\ and Root-node.}
         \label{fig:adv_diff_2}
     \end{subfigure}
     \hfill
     \begin{subfigure}[b]{0.49\textwidth}
         \centering
         \includegraphics[width=\textwidth]{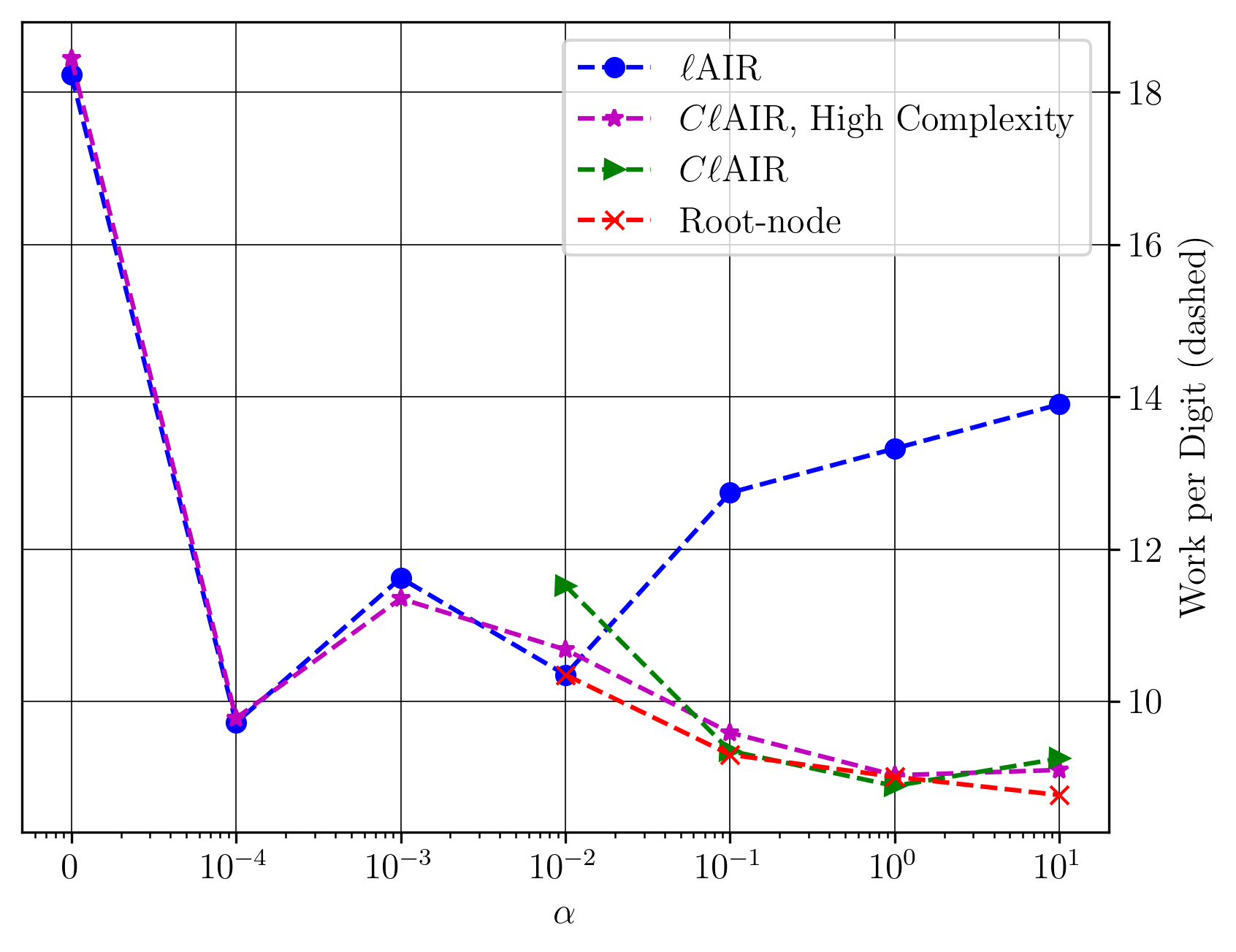}
         \caption{Constant advection with varying diffusion}
         \label{fig:adv_diff_3}
     \end{subfigure}
     \hfill
     \begin{subfigure}[b]{0.49\textwidth}
         \centering
         \includegraphics[width=\textwidth]{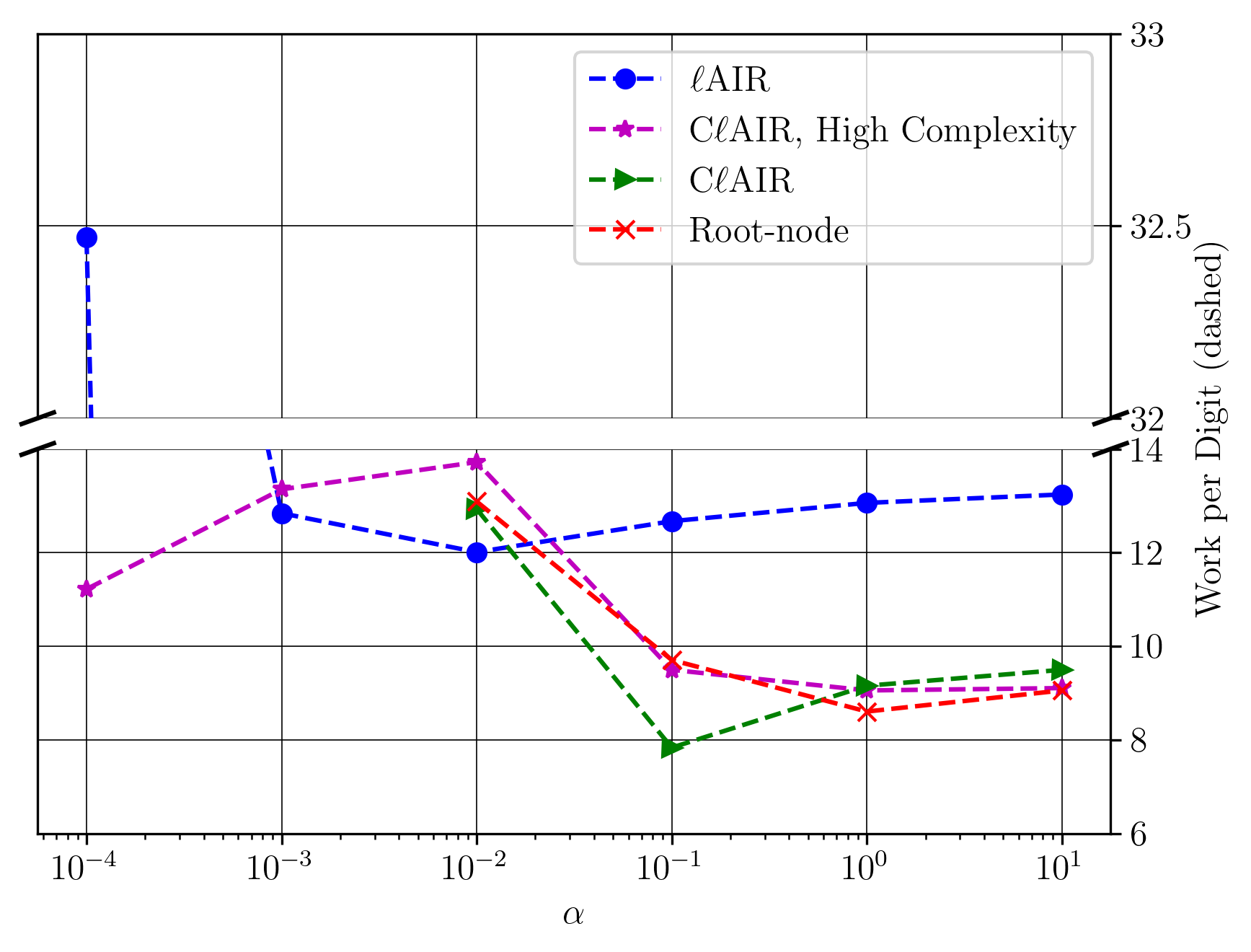}
         \caption{Recirculating advection with varying diffusion}
         \label{fig:adv_diff_4}
     \end{subfigure}
     
        \caption{Advection-diffusion problems, comparison of iterations, work-per-digit of accuracy and operator complexity for $\ell$AIR, \clair\ and Root-node.}
        \label{fig:advection_diffusion_data}
\end{figure}

Figure \ref{fig:adv_diff_1} depicts work-per-digit for the 
2D Constant Direction Advection problem with no diffusion ($\alpha=0$).  The purpose of this
plot is to show that \clair\ enjoys similar convergence to $\ell$AIR for this test problem, where we know that $\ell$AIR works well and is essentially the state-of-the-art \cite{Ben}.

Figures \ref{fig:adv_diff_3} and \ref{fig:adv_diff_4} depict for all three solvers (with two different variants of \clair) the work-per-digit of accuracy over a range of diffusion ($\alpha$) values for the constant direction and recirculating test problems.  The data points for root-node are omitted whenever the solver did not converge within 100 iterations (typically for small $\alpha$ values). We plainly see that \clair\  has the most consistent performance in terms of work-per-digit accuracy across all regimes.

Table \ref{table_op_comp_adv} highlights the operator complexity advantage of \clair\ over $\ell$AIR for the constant advection case (the recirculating free case is similar).  For the more advective cases (smaller $\alpha$), the high complexity variant of \clair, using first-pass only Ruge-St\"{u}ben coarsening, obtains lower
operator complexities of roughly 15--30\%, while for the more diffusive cases (larger $\alpha$), the lower complexity variant of \clair, using aggregation based coarsening, obtains operator complexities roughly 1.6x--2.5x smaller.  In both settings, significant storage savings are achieved. 
\begin{table}[h!]
\begin{center}
\begin{tabular}{ c|c|c|c |c|c|c|c } 
\toprule
\hfill $\alpha$          & 10.0 &  1.0 &  0.1 & 0.01 & 0.001 & 0.0001 & 0.0 \\
\midrule
$\ell$AIR                & 3.19 & 3.19 & 3.20 & 3.19 & 3.70  & 3.19   & 4.28 \\
\clair, High Complexity  & 3.00 & 2.97 & 2.94 & 2.81 & 2.89  & 2.69   & 3.06 \\
\clair                   & 1.28 & 1.29 & 1.29 & 1.96 &  DNC  &  DNC   &  DNC \\
Root-node                & 1.20 & 1.20 & 1.20 & 1.42 &  DNC  &  DNC   &  DNC \\
\bottomrule
\end{tabular}
\end{center}
\caption{Operator complexities of $\ell$AIR and \clair\, for different diffusion $\alpha$ values and the constant advection test problem.  Entries ``DNC" indicate the solver did not converge within 100 iterations.}
\label{table_op_comp_adv}
\end{table}
 
\begin{remark}
For these plots, our goal is to highlight the robustness (lack of need for tuning) for the \clair\ solver.  Thus, we searched for good general parameters for each solver and held them fixed over all tests.  However, better parameters for individual problems and $\ell$AIR do exist.  For instance, if the matrix is \emph{not} block diagonally pre-scaled, then the solver iterations and work for $\ell$AIR in Figure \ref{fig:adv_diff_2} ($\alpha=10$ case) improve and become 7, 8, 8, 9 and 10.35, 10.43, 11.60, 11.81, respectively, but with a prohibitively large operator complexity of $\sim 4.23$.  However with no pre-scaling, the $\ell$AIR then fails to converge within 100 iterations for the small $\alpha$ cases, e.g., Figure \ref{fig:adv_diff_1}.  Another place where tuning is beneficial for $\ell$AIR is the smallest $\alpha$ value for the recirculating advection-diffusion problem ($\alpha = 10^{-4})$, where using first-pass only Ruge-St\"{u}ben coarsening lowers the iterations and work from 88 and 85.34 to 21 and 15.86, respectively.  However, the use of first-pass only Ruge-St\"{u}ben coarsening then substantially degrades iterations and work-per-digit results for $\ell$AIR and other $\alpha$ values.
\end{remark}
\begin{remark}
    For a subset of the advection-diffusion problems considered in this section, Appendix A considers AMG approximation properties of \clair\ and $\ell$AIR, with the results indicating that \clair\ either maintains or improves on the approximation properties of $\ell$AIR.
\end{remark}
In summary, the benefits of \clair\ are as follows.  The $\ell$AIR solver (i) requires more tuning for these problems than \clair, (ii) requires more work-per-digit than \clair\ for more diffusive problems (larger $\alpha$ values), and (iii) has significantly larger operator complexities. 

\section{Conclusion}
In this paper, we developed a new type of reduction-based AMG that is suitable for solving nonsymmetric linear systems coming from the discretization of advection-diffusion PDEs.  By combining techniques from $\ell$AIR that have been effective for advective problems, with energy minimization and root node techniques that are well suited for diffusion problems, we have developed an efficient method for solving advection-diffusion problems in a general setting -- the solver is insensitive to varying contribution of the diffusive part in the PDE.  An important distinction between our proposed solver and existing reduction based methods is that we take a column-wise approach to computing an approximation to ideal restriction and interpolation, which more naturally allows us to incorporate energy minimization and mode constraint techniques into the process.  Our future work focuses on deriving a two-grid convergence theory for the proposed approach applied to nonsymmetric systems and to incorporate the idea of compatible relaxation into the \clair\ coarsening process.

\appendix
\label{sec:appendix_approx}
\section{Classical AMG Weak and Strong Approximation Properties} \subsection{Approximation Properties in the Nonsymmetric Setting}  We consider convergence of $\ell$AIR and \clair\ based on classical multigrid weak and strong approximation properties. Targeting nonsymmetric problems, we consider the generalization of the $A$-norm as $\sqrt{A^{*}A}$ or $\sqrt{AA^{*}}$  \cite{brezina2010towards,BenTheory}. For a nonsingular matrix $A \in \mathbb{C}^{n\times n}$, consider the singular value decomposition (SVD) given by $A=U\Sigma V^{*}$ where the singular values are $0<\sigma_{1}\leq \sigma_{2}\leq...\leq \sigma_{n}$. Then $\sqrt{A^{*}A}=V\Sigma V^{*}=VU^{*}U\Sigma V^{*}= QA$ where $Q=VU^{*}$. In a similar manner, we can also obtain  $\sqrt{AA^{*}}=U\Sigma U^{*}=AQ$. Since $\sqrt{A^{*}A}$ or $\sqrt{AA^{*}}$ are SPD matrices, therefore we can still consider classical AMG approximation properties with respect to the SPD matrices $QA$ and $AQ$ corresponding to the right and left singular vectors, respectively.  Such approximation properties measure, in a sense, how effective the coarse spaces are at capturing the near nullspace of the operator.

For several test problems, we numerically evaluate approximation properties for $\ell$AIR and \clair\ by making use of the generalized fractional approximation property (FAP) \cite{BenTheory}. The FAP of a transfer operator $T \in \mathbb{R}^{n \times n_c}$ is with respect to the SPD matrix $\mathcal{A}$, with powers $\beta, \eta \geq 0$ and constant $K_{T,\beta,\eta}$. Specifically, $T$ is said to have a FAP if for every fine grid vector, $\mathbf{v}$, there exists a coarse grid vector, $\mathbf{v}_{c}$, such that 
\begin{equation}\label{FAP}
||\mathbf{v}-T\mathbf{v}_c||_{\mathcal{A}^{\eta}}^{2}\leq \displaystyle\frac{K_{T,\beta, \eta}}{||\mathcal{A}||^{2\beta-\eta}}\langle\mathcal{A}^{2\beta}\mathbf{v},\mathbf{v}\rangle,
\end{equation}
where $\mathcal{A} = QA$ ($T = P$) or $\mathcal{A} = AQ$ ($T = R^*$). The classical multigrid weak approximation property (WAP) is a FAP($1/2, 0$), that is 
\begin{equation}\label{WAP}
   ||\mathbf{v}-T\mathbf{v}_c||_{2}^{2}\leq \displaystyle\frac{K_{T,1/2, 0}}{||\mathcal{A}||}\langle\mathcal{A}\mathbf{v},\mathbf{v}\rangle.
\end{equation}
Further, the classical multigrid strong approximation property (SAP) is a FAP(1,1), that is 
\begin{equation}\label{SAP}
    ||\mathbf{v}-T\mathbf{v}_c||_{\mathcal{A}}^{2}\leq \displaystyle\frac{K_{T,1, 1}}{||\mathcal{A}||}\langle\mathcal{A}^{2}\mathbf{v},\mathbf{v}\rangle.
\end{equation}
For a given vector $\mathbf{v}$, we compute the approximation constant $K_{T,\beta, \eta}(\mathbf{v})$ with
\begin{equation}\label{minimize}
    K_{T,\beta, \eta}(\mathbf{v})=\displaystyle\frac{||\mathcal{A}||^{2\beta-\eta}}{||\mathbf{v}||^{2}_{\mathcal{A}^{2\beta}}}\,\,\min_{\mathbf{v}_{c}}||\mathbf{v}-T\mathbf{v}_c||_{\mathcal{A}^{\eta}}^{2}.
\end{equation}

\noindent Let $\Pi_{\eta}$ denote the $\mathcal{A}^{\eta}$-orthogonal projection onto the range of $T$, $\Pi_{\eta}= T(T^{*}\mathcal{A}^{\eta}T)^{-1}T^{*}\mathcal{A}^{\eta}  $. 
Substituting into \eqref{minimize} we obtain,

\begin{equation}\label{minimize2}
    K_{T,\beta, \eta}(\mathbf{v})=\displaystyle\frac{||\mathcal{A}||^{2\beta-\eta}}{||\mathbf{v}||^{2}_{\mathcal{A}^{2\beta}}}\,\,||(I-\Pi_{\eta})\mathbf{v}||_{\mathcal{A}^{\eta}}^{2}.
\end{equation}

\noindent
To compute the approximation constant $K_{\max}$ that holds for all fine grid vectors $\mathbf{v}$, we take maximum of the above expression over all $\mathbf{v}$, which leads to
\begin{equation}    
  K_{\max}=\max_{\mathbf{v}\neq 0}\,\,K_{T,\beta, \eta}(\mathbf{v})={||\mathcal{A}||^{2\beta-\eta}}\,\,||\mathcal{A}^{\eta/2} (I-\Pi_{\eta}) \mathcal{A}^{-\beta}||_{2}^{2}.
\end{equation}

\subsection{Numerical Tests} To complement the numerical results in the main text, we now measure classical AMG approximation property constants for both $\ell$AIR and \clair\ from the purely advective to diffusion dominated case. Specifically, we consider the 2D constant direction advection problem \eqref{eqn:advdiff1}--\eqref{eqn:advdiff2} for various diffusion coefficients $\alpha$. For numerical tests, a $32\times 32$ size spatial domain is considered, resulting in $1024$ total DOFs. In the following numerical tests, we consider the approximation properties of the restriction operators from both $\ell$AIR and \clair.

\begin{figure}[h!]
     \centering
     \begin{subfigure}[b]{0.49\textwidth}
         \centering
         \includegraphics[width=\textwidth]{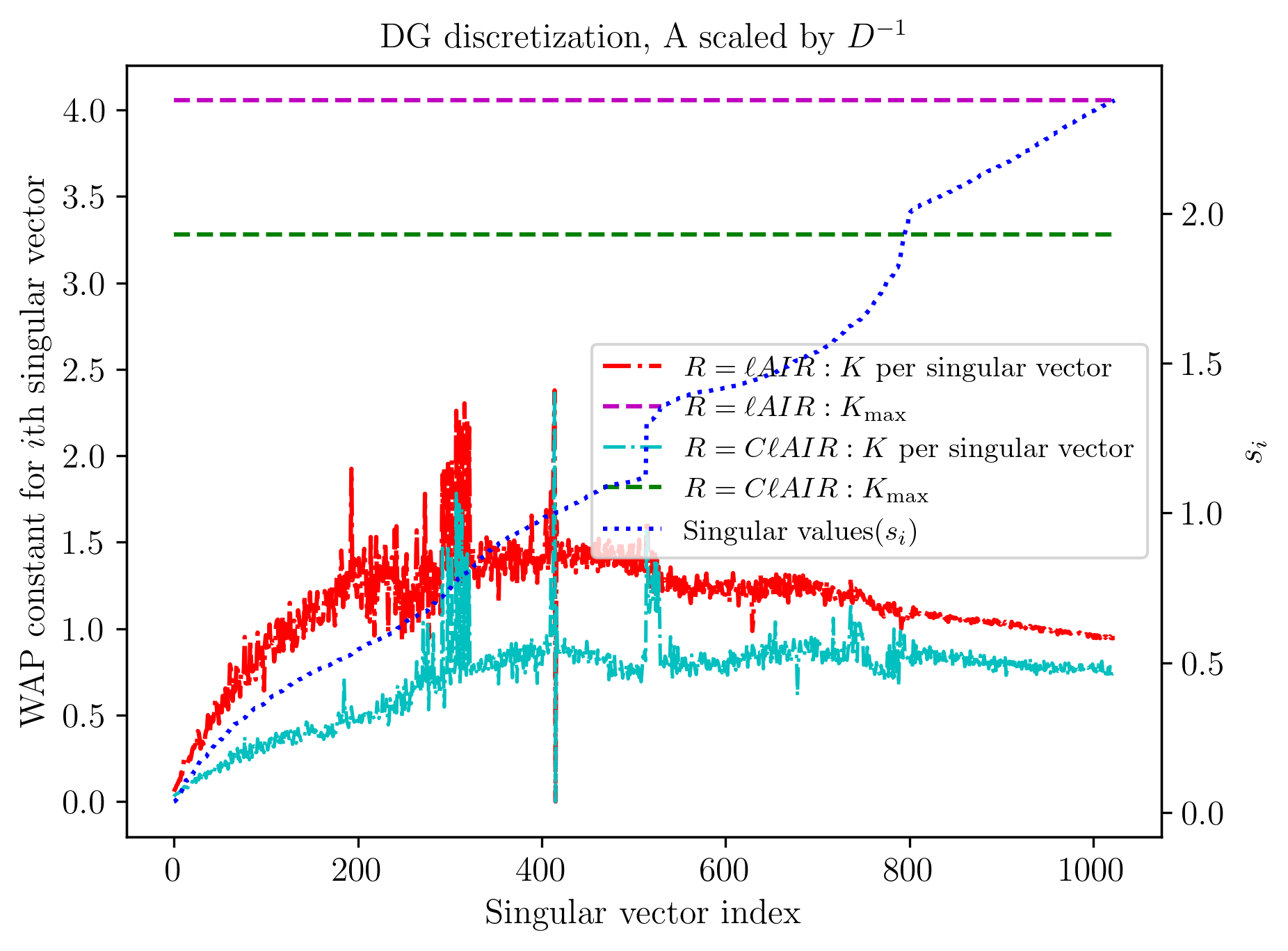}
         \caption{WAP for Constant Advection with zero diffusion ($\alpha=0$) where constraint vector $B = \mathbf{1}$ is presmoothed with 5 iterations of CFF-weighted-Jacobi.}
         \label{fig:air_approx_0}
     \end{subfigure}
     \hfill
     \begin{subfigure}[b]{0.49\textwidth}
         \centering
         \includegraphics[width=\textwidth]{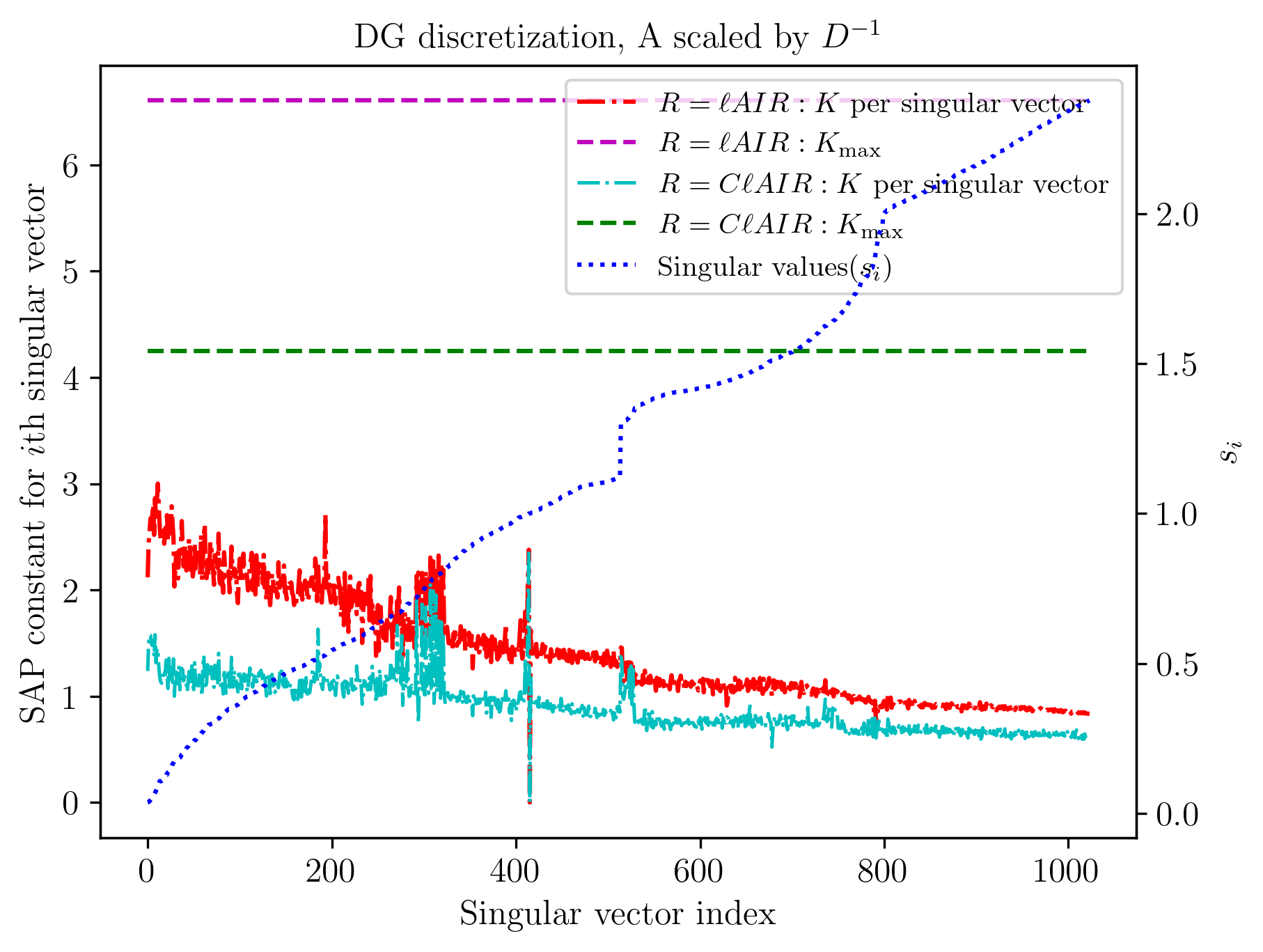}
         \caption{SAP for Constant Advection with zero diffusion ($\alpha=0$) where constraint vector $B = \mathbf{1}$ is presmoothed with 5 iterations of CFF-weighted-Jacobi.}
         \label{fig:air_approx_1}
     \end{subfigure}

        \caption{WAP and SAP constants for the restriction operators of $\ell$AIR and \clair\ for the constant advection with 
 zero diffusion ($\alpha=0$) problem. Five (5) iterations of CFF-weighted-Jacobi relaxation has been used to improve the mode constraint vector $B=\mathbf{1}$ in \clair. Singular values are shown in the dotted blue line and are associated with the right vertical axis, and the dot-dashed lines show the approximation constant for each of the left singular vectors of $A$. Horizontal dashed lines show the approximation constant for $\ell$AIR and \clair\ that holds for all vectors.}  
        \label{fig:approx_fig_air_zero_diff}
\end{figure}

Figure \ref{fig:approx_fig_air_zero_diff} shows the WAP (FAP($1/2, 0$)) and SAP (FAP($1,1$)) approximation constants for the left singular vectors of $A$. Here, 5 iterations of CFF-weighted-Jacobi relaxation has been used to improve the mode constraint vector $B=\mathbf{1}$. From Figure \ref{fig:air_approx_0} and \ref{fig:air_approx_1}, we see that \clair\ demonstrates moderately smaller (i.e., better) WAP and SAP constants than $\ell$AIR. Since $\ell$AIR is primarily designed for such purely advective problems these results indicate that \clair\ is also effective for advection problems, with similar convergence properties. This is indeed what we see in Figures \ref{fig:adv_diff_1} and \ref{fig:adv_diff_3}. This motivates us to further study approximation properties of \clair\ for  diffusion dominated problems.

\begin{figure}[h!]
     \centering
     \begin{subfigure}[b]{0.49\textwidth}
         \centering
         \includegraphics[width=\textwidth]{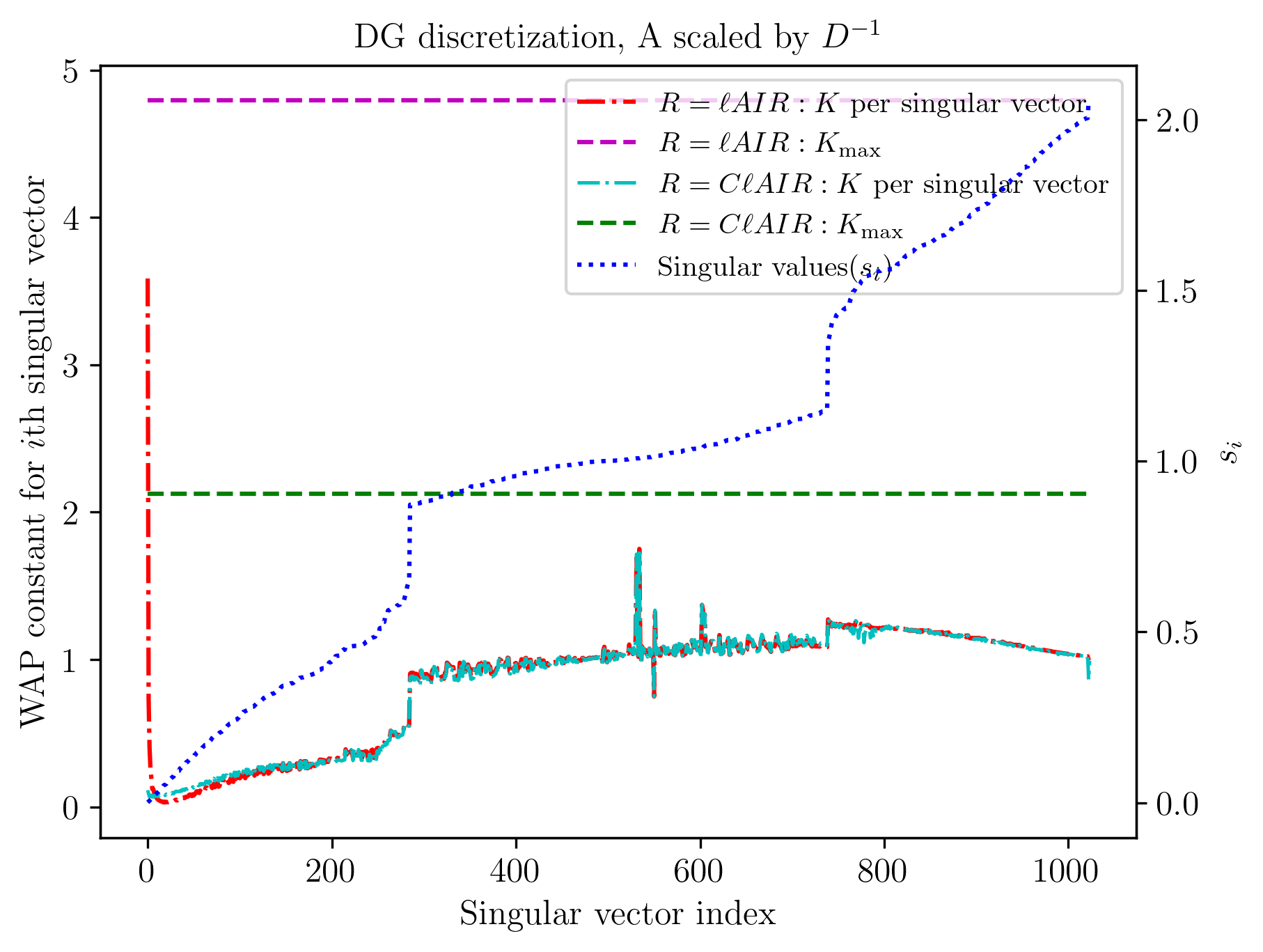}
         \caption{WAP for Constant Advection with added diffusion ($\alpha=10$) where constraint vector $B = \mathbf{1}$ is presmoothed with 5 iterations of CFF-weighted-Jacobi.}
         \label{fig:air_approx_2}
     \end{subfigure}
     \hfill
     \begin{subfigure}[b]{0.49\textwidth}
         \centering
         \includegraphics[width=\textwidth]{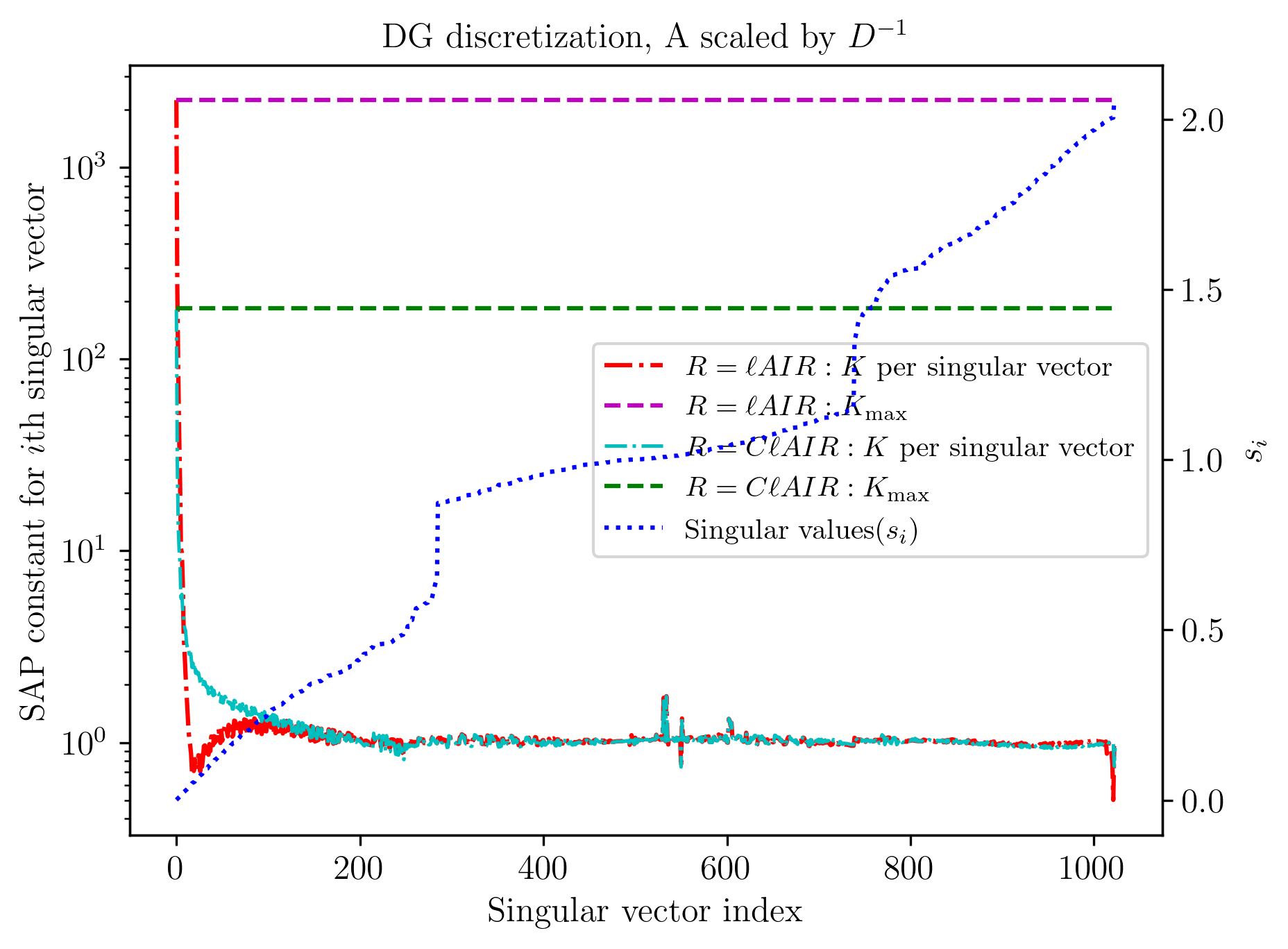}
         \caption{SAP for Constant Advection with added diffusion ($\alpha=10$) where constraint vector $B = \mathbf{1}$ is presmoothed with 5 iterations of CFF-weighted-Jacobi.}
         \label{fig:air_approx_3}
     \end{subfigure}
     \begin{subfigure}[b]{0.49\textwidth}
         \centering
         \includegraphics[width=\textwidth]{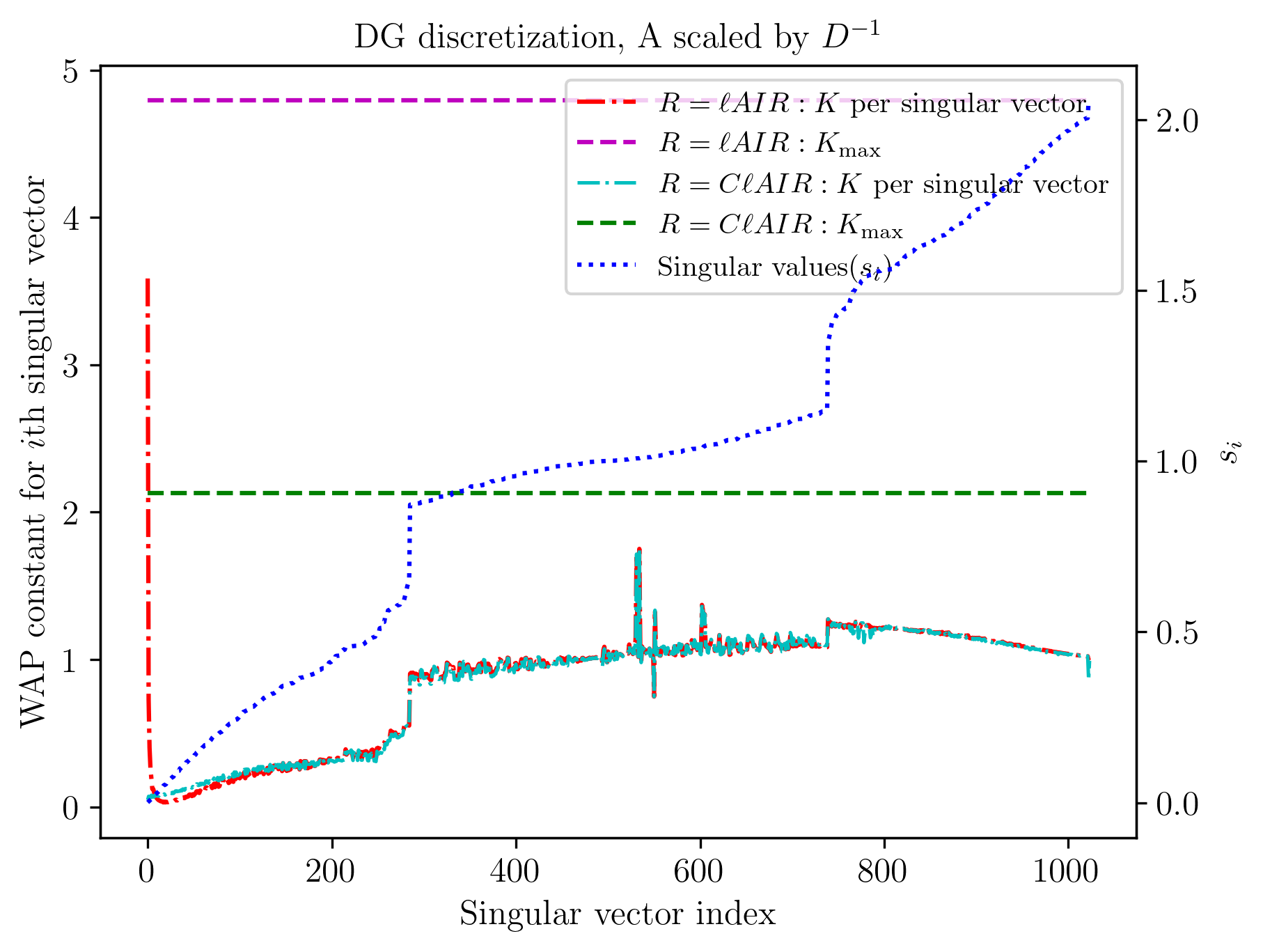}
         \caption{WAP for Constant Advection with added diffusion ($\alpha=10$) where constraint vector $B = \mathbf{1}$ is presmoothed with 25 iterations of CFF-weighted-Jacobi.}
         \label{fig:air_approx_4}
     \end{subfigure}
        \hfill
     \begin{subfigure}[b]{0.49\textwidth}
         \centering
         \includegraphics[width=\textwidth]{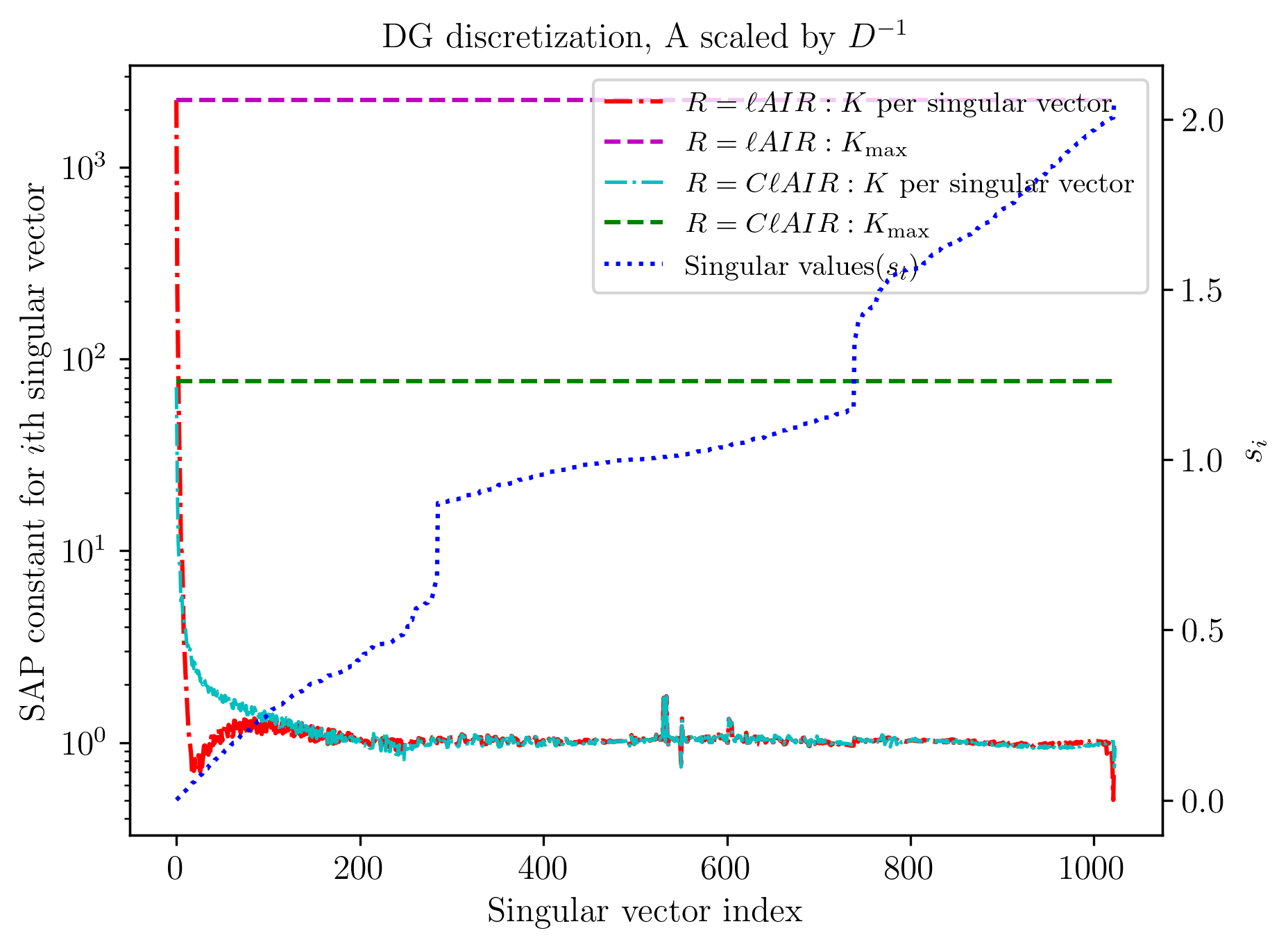}
         \caption{SAP for Constant Advection with added diffusion ($\alpha=10$) where constraint vector $B = \mathbf{1}$ is presmoothed with 25 iterations of CFF-weighted-Jacobi.}
         \label{fig:air_approx_5}
     \end{subfigure}
     
        \caption{WAP and SAP constants for the restriction operators of $\ell$AIR and \clair\ for the constant advection with 
 added diffusion ($\alpha=10$) problem. Different number of iterations (5 and then 25) of CFF-weighted-Jacobi relaxation has been used to improve the mode constraint vector $B=\mathbf{1}$ in \clair. Singular values are shown in the dotted blue line and are associated with the right vertical axis, and the dot-dashed lines show the approximation constant for each of the left singular vectors of $A$. Horizontal dashed lines show the approximation constant for $\ell$AIR and \clair\ that holds for all vectors.}  
 \label{fig:approx_fig_air_10_diff}
 \end{figure}

Next, we consider a diffusion dominated case by setting the  diffusion coefficient to be $\alpha=10$. First, as previously, we use 5 iterations of a CFF-weighted-Jacobi smoother to improve the mode constraint vector $B=\mathbf{1}$. Carefully investigating the results in Figure \ref{fig:air_approx_2}, we observe that WAP constant for the new method \clair\ is almost a factor of two smaller than that for $\ell$AIR. Now if we compare the SAP constants in Figure \ref{fig:air_approx_3}, $K_{\max}$ is smaller for \clair\ than for $\ell$AIR, although the constants are very large for both solvers (184 for \clair\ and 2246 for $\ell$AIR). Somewhat surprisingly, while the SAP constant for \clair\ in Figure \ref{fig:air_approx_3} is quite large, we find that  the iteration counts of the solver are mesh independent (as we have seen in Figure \ref{fig:adv_diff_2} for the $\alpha=10$ case). In this case, we suspect that the Krylov method is able to account for what the solver is lacking. We note that the significantly larger SAP constants for $\ell$AIR are consistent with the poorer scalability of $\ell$AIR seen in Figure \ref{fig:adv_diff_2} for the $\alpha=10$ case. A key advantage of \clair\ is that it allows the flexibility to improve the approximation properties of the restriction operator by employing a suitable relaxation scheme for smoothing the mode constraint vector (such as weighted CFF-Jacobi). This type of flexibility is not available in $\ell$AIR. Therefore, further improvement of \clair's SAP constant can be obtained by increasing the number of iterations applied to $B$. In our experiments, we find that if the mode constraint vector $B = \mathbf{1}$ is presmoothed with 25 iterations of CFF-weighted-Jacobi (instead of 5 iterations), the SAP approximation constant for \clair\ decreases further to $K_{\max} = 77$ (Figure \ref{fig:air_approx_5}).  However, we do not find that these extra iterations lead to a meaningful improvement in practice for \clair\ convergence on the tested problems.

\subsection*{Acknowledgments}
Los Alamos National Laboratory Report LA-UR-23-27101.

  \bibliography{refs}
  \bibliographystyle{plain}

\end{document}